\documentclass[11pt]{amsart}
\usepackage{amsmath}
\usepackage{amsfonts}
\usepackage{latexsym}
\usepackage{amssymb}
\usepackage{enumerate}
\usepackage[dvips]{graphics}

\newcommand{\Z}{{\mathbf Z}}
\newcommand{\Q}{{\mathbf Q}}
\newcommand{\R}{{\mathbf R}}
\newcommand{\C}{{\mathbf C}}

\newcommand{\kk}{{\mathbf k}}
\newcommand{\cat}{{\rm {cat }}}

\newcommand{\Hom}{{\rm Hom}}
\newcommand{\im}{\rm im}

\newcommand{\ind}{{\rm {ind}}}
\newcommand{\Int}{{\rm {Int}}}
\newcommand{\Supp}{{\rm {Supp}}}
\newcommand{\ccat}{{\rm {ccat}}}
\newcommand{\cl}{{\rm {cl}}}

\newcommand{\rk}{{\rm rk}}
\newcommand{\comment}[1]{}

\def\V{\mathcal V}

\def\ker{{\rm Ker }}
\def\im{{\rm Im }}

\def\R{\mathbf R}

\newtheorem{theorem}{Theorem}
\newtheorem{question}{Question}
\newtheorem{proposition}{Proposition}[section]
\newtheorem{lemma}[proposition]{Lemma}
\newtheorem{example}[proposition]{Example}
\newtheorem{corollary}[proposition]{Corollary}
\newtheorem{definition}[proposition]{Definition}
\newtheorem{remark}[proposition]{Remark}

\begin{document}
\title{Cohomological estimates for {\Large {$\cat(X,\xi)$}}}
\author[M. ~Farber and D. ~Sch\"utz]{M.~Farber and D. ~Sch\"utz}
\address{Department of Mathematics, University of Durham, Durham DH1 3LE, UK}
\email{Michael.Farber@durham.ac.uk}

\email{Dirk.Schuetz@durham.ac.uk}



\date{}


\begin{abstract} This paper studies the homotopy invariant $\cat(X,\xi)$ introduced in
\cite{farbe2}. Given a finite cell-complex $X$, we study the function
$\xi\mapsto \cat(X,\xi)$ where $\xi$ varies in the cohomology space
$H^1(X;\R)$. Note that $\cat(X,\xi)$ turns into the classical Lusternik -
Schnirelmann category $\cat(X)$ in the case $\xi=0$. Interest in $\cat(X,\xi)$
is based on its applications in dynamics where it enters estimates of
complexity of the chain recurrent set of a flow admitting Lyapunov closed
1-forms, see \cite{farbe2}, \cite{farbe3}.

In this paper we significantly improve earlier cohomological lower bounds for
$\cat(X,\xi)$ suggested in \cite{farbe2}, \cite{farbe3}. The advantages of the
current results (see Theorems \ref{thm1}, \ref{thm11} and \ref{perfect} below)
are twofold: firstly, we allow cohomology classes $\xi$ or arbitrary rank
(while in \cite{farbe2} the case of rank one classes was studied), and
secondly, the theorems of the present paper are based on a different principle
and give slightly better estimates even in the case of rank one classes. We
introduce in this paper a new controlled version of $\cat(X,\xi)$ and find
upper bounds for it (Theorems \ref{thm3} and \ref{improved}). We apply these
upper and lower bounds in a number of specific examples where we explicitly
compute $\cat(X,\xi)$ as a function of the cohomology class $\xi\in H^1(X;\R)$.
\end{abstract}

\maketitle
\section{Introduction}

This paper studies the homotopy invariant $\cat(X,\xi)$ originally introduced
in \cite{farbe2} and investigated further in \cite{farbe3}. It is a
generalization of the classical Lusternik - Schnirelman category which appeared
in connection with the problem of estimating the number of zeros of closed
one-forms lying in a prescribed nonzero cohomology class. Assuming that the
zeros are all Morse type the answer is provided by the Novikov theory
\cite{noviko}. If the zeros are arbitrary, there always exists a closed 1-form
with at most one zero \cite{farbe2}, \cite{FS1} and therefore one cannot expect
to get meaningful lower bounds on the number of zeros in homotopic terms.
However, such lower bounds exist if the gradient flow of the closed one form
has no homoclinic cycles \cite{farbe2}, \cite{farbe3}. More generally, there is
always an interesting relation, governed by $\cat(X,\xi)$, between the topology
of the chain recurrent set of the flow admitting a Lyapunov 1-form and the
topology of the underlying manifold, see \cite{FK}, \cite{Lat}.

Our goal in this paper is to find a cohomological lower bound for $\cat(X,\xi)$
generalizing the classical cup-length estimate for the usual category. Earlier
cohomological lower bounds for $\cat(X,\xi)$ (see \cite{farbe2}, Theorems 6.1
and 6.4) were obtained only in the case of rank one cohomology classes. The
advantages of the current results (see Theorems \ref{thm1}, \ref{thm11} and
\ref{perfect} below) are twofold: firstly, we allow cohomology classes $\xi$ or
arbitrary rank and secondly, the theorems of the present paper are based on a
different principle and give slightly better estimates even in the case of rank
one classes.

The main result of the paper is Theorem \ref{perfect} which we will now
describe. Let $(X,\xi)$ be a pair consisting of a finite cell-complex and a
cohomology class $\xi\in H^1(X;\R)$. We denote by $\ker(\xi)\subset H_1(X;\Z)$
the subgroup consisting of homology classes $z$ such that $\langle \xi,
z\rangle =0$. The factor group $H=H_1(X;\Z)/\ker (\xi)$ is a free abelian group
of rank $r=\rk(\xi)$. We consider the set $\V_\xi$ of all complex flat line
bundles $L$ over $X$ such that the monodromy of $L$ is trivial along any
homology class in $\ker(\xi)$. Clearly, $\V_\xi =\Hom(H;\C^\ast) = (\C^\ast)^r$
has the structure of an algebraic variety. A bundle $L\in \V_\xi$ is called
transcendental if the monodromy homomorphism ${\rm {Mon}}_L: \Z[H]\to \C$ is
injective. In the case when $\rk H=1$ line bundles correspond to complex
numbers (since $\V_\xi=\C^\ast$) and transcendental bundles correspond to
numbers which are transcendental in the usual sense\footnote{Numbers which are
not roots of polynomial equations with integral coefficients}. One of our main
results is given by the following theorem:

\begin{theorem}
Suppose that $L\in \V_\xi$ is transcendental and there exist cohomology classes
$v_0\in H^\ast(X;L)$ and $v_i\in H^{d_i}(X;\C)$ where $i=1, \dots, k$. If
$d_i>0$ for $i=1, \dots, k$ and the cup-product $v_0\cup v_1\cup\dots\cup
v_k\in H^\ast(X;L)$ is nonzero then $\cat(X,\xi)>k.$
\end{theorem}

Given a finite cell-complex $X$, we examine the function $\xi\mapsto
\cat(X,\xi)$ where $\xi$ varies in the cohomology space $H^1(X;\R)$. We compute
it explicitly in some examples and investigate its behavioral patterns.

Another new result worth mentioning is the introduction of a new controlled
version of $\cat(X,\xi)$ which we denote $\ccat(X,\xi)$. It is a modification
of the invariant $\cat(X,\xi)$ which coincides with $\cat(X,\xi)$ in all
examples known to us. We find new upper bounds for $\ccat(X,\xi)$ (Theorems
\ref{thm3} and \ref{improved}) and a product inequality generalizing the
classical product inequality for the usual category.

There are interesting connections between the invariant $\cat(X,\xi)$ and its
various modifications and the Bieri-Neumann-Strebel invariant of discrete
groups \cite{BNS}, see also \cite{BG}. We discuss some of these relations in \S
\ref{bieri}. We show that the Bieri-Neumann-Strebel invariant allows us to
improve the upper bound for $\cat(X,\xi)$ in the case of manifolds. In
particular we prove the following theorem (which is combination of Theorems
\ref{thm3} and \ref{ahh} in the text):

\begin{theorem}\label{ahhh}
Let $M$ be a closed connected smooth manifold. Then for any nonzero cohomology
class $\xi\in H^1(M;\R)$ one has $\cat(M,\xi) \leq n-1$ where $n=\dim M$.
Moreover, if $n\geq 5$ and for some nonzero $\xi\in H^1(M;\R)$ we have
$\cat(M,\xi) = n-1$ then the fundamental group $\pi_1(M)$ contains a
non-abelian free subgroup.
\end{theorem}

In the last section of the paper we discuss several open problems.
\section{Abel - Jacobi maps and neighborhoods of infinity}\label{sec:aj}

For convenience of the reader we give in sections \ref{sec:aj} and
\ref{sec:movable} a brief summary of our result \cite{FS} (see Theorem
\ref{field} below) which will play a crucial role in this paper.

Let $X$ be a connected finite cell complex and $p:\tilde X\to X$ a regular
covering having a free abelian group of covering transformations $H\simeq
\Z^r$. Denote $ H_\R=H\otimes \R; $ it is a vector space of dimension $r$
containing $H$ as a lattice.
\begin{proposition}
There exists a canonical Abel - Jacobi map
\begin{eqnarray} A: \tilde X\to H_\R
\end{eqnarray}
having the following properties:
\begin{enumerate}
\item[(a)] $A$ is $H$-equivariant; here $H$ acts on $\tilde X$ by covering
transformations and it acts on $H_\R$ by translations.

 \item[(b)] $A$ is proper (i.e. the preimage of
a compact subset of $H_\R$ is compact).

\item[(c)] $A$ is determined uniquely up to replacing it by a map $A': \tilde
X\to H_\R$ of the form $ A' = A+ F\circ p $ where $F: X\to H_\R$ is a
continuous map.
\end{enumerate}\end{proposition}

This fact is well-known; see for example \cite{BNS}, \cite{FS}.

Let $\xi\in H^1(X;\R)$ be a cohomology class with the property
$$p^\ast(\xi)=0\in H^1(\tilde X;\R).$$ Such a class $\xi$ can be viewed either as a
homomorphism $\xi: H\to \R$ or as a linear functional $\xi_\R: H_\R\to \R$.

\begin{definition}
A subset $N\subset \tilde X$ is called a neighborhood of infinity in $\tilde X$
with respect to the cohomology class $\xi$ if $N$ contains the set
\begin{eqnarray}
\{x\in \tilde X; \, \, \xi_\R(A(x))>c\}\, \subset \, N,
\end{eqnarray}
for some real $c\in \R$. Here $A: \tilde X\to H_\R$ is an Abel - Jacobi map for
the covering $p:\tilde X\to X$.
\end{definition}

This term was introduced in \cite{FS}. However, the object itself is certainly
not new. It appeared earlier (with no name) in \cite{BNS} and in other papers.

Neighborhoods of infinity $N\subset \tilde X$ with respect to $\xi$ can be
characterized as follows. Let $\gamma: S^1\to X$ be a continuous map such that
evaluation of class $\xi$ on the homology class $[\gamma]\in H_1(X;\Z)$ is
positive, $\langle \xi, [\gamma]\rangle >0$. Consider a lift $\tilde \gamma :
\R\to \tilde X$ of the composition $\R\stackrel \exp\to S^1\stackrel \gamma\to
\tilde X$, see the commutative diagram:
\begin{eqnarray}
\begin{array}{rcl}
\R & \stackrel {\tilde \gamma}\to & \tilde X\\
{\rm {exp}}\downarrow& & \downarrow p\\
S^1&\stackrel \gamma \to & X.
\end{array}
\end{eqnarray}
Then for all sufficiently large $t\in \R$ the point $\tilde \gamma(t)$ lies in
the neighborhood $N\subset \tilde X$.

\section{Homology classes movable to infinity}\label{sec:movable}

Let $G$ be an abelian group (the coefficient system). We mainly have in mind
the cases $G=\Z$ of $G=\kk$ is a field.

\begin{definition}
A homology class $z\in H_i(\tilde X;G)$ is said to be movable to infinity of
$\tilde X$ with respect to a nonzero cohomology class $\xi\in H^1(X;\R)$,
$p^\ast(\xi)=0$, if in any neighborhood $N$ of infinity with respect to $\xi$
there exists a (singular) cycle with coefficients in $G$ representing $z$.
\end{definition}

Equivalently, a homology class $z\in H_i(\tilde X;G)$ is movable to infinity
with respect to $\xi\in H^1(X;\R)$ if $z$ lies in the intersection
\begin{eqnarray}
\bigcap_N  \im\left[ H_i(N;G)\to H_i(\tilde X;G)\right]
\end{eqnarray}
where $N$ runs over all neighborhoods of infinity in $\tilde X$ with respect to
$\xi$. This can also be expressed by saying that $z$ lies in the kernel of the
natural homomorphism
\begin{eqnarray}\label{kernel}
H_i(\tilde X;G) \to \lim_\leftarrow H_i(\tilde X,N;G)
\end{eqnarray}
where in the inverse limit $N$ runs over all neighborhoods of infinity in
$\tilde X$ with respect to $\xi$.


The following theorem proven in \cite{FS} gives an explicit description of all
movable homology classes in the case when $G=\kk$ is a field. It generalizes
the result of \cite{farbe2}, \S 5 treating the simplest case of infinite cyclic
covers $q: \tilde X\to X$.

\begin{theorem}\label{field} Let $X$ be a finite cell complex and $q:\tilde X\to X$
be a regular covering having a free abelian group of covering transformations
$H\simeq\Z^r$. Let $\xi\in H^1(X;\R)$ be a nonzero cohomology class of rank $r$
satisfying $q^\ast(\xi)=0$. The following properties (A), (B), (C) of a nonzero
homology class $z\in H_i(\tilde X;\kk)$ (where $\kk$ is a field) are
equivalent:

\begin{enumerate}
\item[(A)] $z$ is movable to infinity with respect to $\xi$.

\item[(B)] Any singular cycle $c$ in $\tilde X$ realizing the class $z$ bounds
an infinite singular chain $c'$ in $\tilde X$ containing only finitely many
simplices lying outside every neighborhood of infinity $N\subset \tilde X$ with
respect to $\xi$.

\item[(C)] There exists a nonzero element $x\in \kk[H]$ such that $x\cdot z=0$.
\end{enumerate}
\end{theorem}

\section{Obstructions to movability to infinity of homology classes}

In this section we show that {\it generic} flat line bundles can be used to
detect movability to infinity of homology classes. This observation is used
later in this paper when we study cohomological lower bounds for $\cat(X,\xi)$,
see \S \ref{sec4}.

Let $X$ be a finite polyhedron and $\xi \in H^1(X;\R)$. Denote by $\ker (\xi)$
the kernel of the homomorphism $H_1(X;\Z)\to \R$ given by evaluation on $\xi$.
Then $H=H_1(X;\Z)/\ker(\xi)$ is a free abelian group of finite rank $r$ where
$r$ denotes the rank of $\xi$. Consider the cover $p: \tilde X\to X$
corresponding to $\ker(\xi)$. It has $H$ as the group of covering translations.

Let $\V_\xi= (\C^\ast)^r= \Hom(H, \C^\ast)$ denote the variety of all complex
flat line bundles $L$ over $X$ such that the induced flat line bundle $p^\ast
L$ on $\tilde X$ is trivial. If $t_1, \dots, t_r\in H$ is a basis, then the
monodromy of $L\in \mathcal V_\xi$ along $t_i$ is a nonzero complex number
$x_i\in \C^\ast$ and the numbers $x_1, \dots, x_r\in \C^\ast$ form a coordinate
system on $\mathcal V_\xi$. Given a flat line bundle $L\in \V_\xi$ the
monodromy representation of $L$ is the ring homomorphism
\begin{eqnarray}\label{monodromy}{\rm
Mon}_L: \C[H]\to \C\end{eqnarray}
sending each $t_i\in H$ to $x_i\in \C^\ast$.
The dual bundle $L^\ast\in \V_\xi$ is such that $L\otimes L^\ast$ is trivial;
if $x_1, \dots, x_r\in \C^\ast$ are coordinates of $L$ then $x_1^{-1}, \dots,
x_r^{-1}\in \C^\ast$ are coordinates of $L^\ast$.

Any nontrivial element $p\in \C[H]$ lying in the kernel of ${\rm Mon}_L$ can be
viewed as a (Laurent) polynomial equation between the variables $x_1, \dots,
x_r$. Alternatively, we will consider algebraic subvarieties $V\subset \V_\xi$.
Any such $V$ is the set of all solutions of a system of equations of the form
$$q_i(x_1, \dots, x_r, x_1^{-1}, \dots, x_r^{-1})=0, \quad i=1, \dots, m$$
where $q_i$ is a Laurent polynomial with complex coefficients $$p_i\in \C[x_1,
\dots, x_r, x_1^{-1}, \dots, x_r^{-1}].$$ This is equivalent to fixing an ideal
$J\subset \Q[H]$ and studying the set of all flat line bundles $L\in \V_\xi$
such that ${\rm Mon}_L(J)=0.$

We are chiefly interested in the subset $A\subset \V_\xi$ consisting of all
line flat bundles $L\in \V_\xi$ satisfying the following property: {\it for any
homology class $z\in H_\ast(\tilde X; \C)$ which is movable to infinity with
respect to $\xi$ one has $$p_\ast(z) =0 \in H_\ast(X;L)$$ where
\begin{eqnarray}\label{proj}
p_\ast: H_\ast(\tilde X;\C)\to H_\ast(X;L)
\end{eqnarray}
is the map induced by the projection $p$}. Homomorphism (\ref{proj}) exists
because the induced flat line bundle $p^\ast(L)$ over $\tilde X$ is trivial,
$p^\ast(L)=\C$.

\begin{lemma}\label{one} The complement $V = \mathcal V_\xi - A$ is contained in
a proper algebraic subvariety of $\V_\xi$.
\end{lemma}

\begin{proof}
Let $T_q\subset H_q(\tilde X;\C)$ denote the subset of $\C[H]$-torsion homology
classes. According to Theorem \ref{field}, $T_q$ is exactly the set of homology
classes in $H_q(\tilde X;\C)$ which are movable to infinity with respect to
$\xi$.

One has $V=\cup V_q$ where $L\in V_q$ if and only if the composition
$$T_q\to H_q(\tilde X;\C) \to H_q(X;L)$$ is
nontrivial.

Let $L_0$ denote the fiber of $L$ over the base point $x_0\in X$.
 It is a one-dimensional complex vector space and the
monodromy representation of $L$ determines a structure of a right
$\Lambda$-module on $L_0$; here $\Lambda=\C[H]$. Homomorphism (\ref{proj})
equals the composition $$H_q(\tilde X;\C) \to H_q(\tilde
X;\C)\otimes_{\Lambda}L_0 \to H_q(X;L).$$ Hence we find that $L\in V_q$ implies
$T_q\otimes_{\Lambda}L_0\not=0$.

 Let $\dots \to
F_1\stackrel d\to F_0\to T_q\to 0$ be a finitely generated free
$\C[H]$-resolution of $T_q$. We obtain: $L\in V_q$ implies that the linear map
\begin{eqnarray}\label{tensor}
F_1\otimes_{\Lambda}L_0\stackrel d\to F_0\otimes_{\Lambda}L_0 \end{eqnarray} is
not onto.

Let $M$ be the square matrix with entries in $\Lambda$ representing $d: F_1\to
F_0$. Then $V_q$ lies in the subvariety $V'_q$ of $\mathcal V_\xi$ given by
equating to zero all minors of $M$ of size $n_0=\rk F_0$. Each of these minors
is a polynomial in variables $x_1, \dots, x_r$ with complex coefficients.

Now we show that $V'=\cup V'_q\subset \mathcal V_\xi$ is proper. Since $T_q$ is
 finitely generated and torsion there exists a nonzero Laurent polynomial
 $\Delta_q \in \C[H]$ such that the multiplication by $\Delta_q$ annihilates
 $T_q$. Looking at the commutative diagram
\begin{eqnarray*}
\begin{array}{crclclcc}
\to & F_1 & \stackrel d\to & F_0& \to & T_q&\to 0\\ \\
&\Delta_q \downarrow & g\swarrow & \downarrow \Delta_q & &\downarrow 0 &  &\\ \\
\to & F_1 & \underset{d}\to & F_0& \to & T_q&\to 0
\end{array}
\end{eqnarray*}
we find that there exists a $\Lambda$-morphism $g: F_0\to F_1$ such that
 $d\circ g$ coincides with the multiplication by $\Delta_q$. It is easy to see that
 $L\notin V'_q$ assuming that
 $\Delta_q(x_1, \dots, x_r) \not=0\in \C$. Hence $V'$ is a proper subset of $\mathcal V_\xi$.
\end{proof}

\begin{definition}\label{defsupp}
The set $\V_\xi-A$ is called the support of $(X,\xi)$. It is denoted by
$\Supp(X,\xi)$.
\end{definition}

As follows from the proof of the previous lemma, a flat line bundle $L\in
\V_\xi$ lies in $\Supp(X,\xi)$ if and only if there exists a $\C[H]$-torsion
homology class $z\in H_q(\tilde X,\C)$ with $p_\ast(z)\not=0\in H_q(X;L)$. Here
$p: \tilde X\to X$ is the free abelian covering corresponding $\ker (\xi)$. By
Lemma \ref{one}, $\Supp(X,\xi)\subset \V_\xi$ is contained in a proper
algebraic subvariety.

\begin{example} {\rm In the case when $\xi$ has rank one, $r=1$, the variety $\V_\xi$
coincides with $\C^\ast=\C -\{0\}$ and the support $\Supp(X,\xi)\subset
\C^\ast$ is a finite set. In this special case Definition \ref{defsupp}
coincides with the definition given in \cite{farbe3}, page 182.}
\end{example}

\begin{remark}\label{rmk1} {\rm Assume that the class $\xi\not=0$ in nonzero and the polyhedron
$X$ is connected. Then the trivial line bundle $\C\in \V_\xi$ belongs to the
support $\Supp(X,\xi)$. Indeed, the class $1\in H_0(\tilde X;\C)$ is torsion
and $p_\ast(1)=1\in H_0(X;\C)$.}
\end{remark}

As a corollary of Lemma \ref{one} we obtain the following statement:

\begin{theorem}\label{cor4}
Let $L$ be a flat line bundle $L\in \mathcal V_\xi$ which does not lie in the
support $\Supp(X,\xi)$. Suppose that $\langle v, p_\ast(z)\rangle \not=0$ where
$v\in H^q(X;L)$, $z\in H_q(\tilde X;\C)$ and $p_\ast(z)\in H_q(X;L^\ast)$. Here
$L^\ast\in \V_\xi$ denotes the flat line bundle dual to $L$. Then $z$ is not
movable to infinity of $\tilde X$ with respect to $\xi$.
\end{theorem}

This Theorem will be used in the next section in the proof of Theorem
\ref{thm1}.

\section{A cohomological lower bound for $\cat(X,\xi)$.}\label{sec4}

Let $X$ be a finite simplicial polyhedron and let $\xi\in H^1(X;\R)$ be a
cohomology class. Our goal is to estimate from below the number $\cat(X,\xi)$
introduced in \cite{farbe2}. For convenience of the reader we recall the
definition of $\cat(X,\xi)$.

Let $\omega$ be a closed 1-form on $X$ representing $\xi$; we use the formalism
of closed 1-forms on topological spaces suggested in \cite{farbe2}; see also
\cite{farbe3}.

\begin{definition}\label{movable}
A subset $A\subset X$ is $N$-movable with respect to $\omega$ (where $N\in \Z$
is an integer) if there exists a continuous homotopy $h_t:A\to X,$ $t\in
[0,1],$ such that $h_0:A\to X$ is the inclusion and for any point $x\in A$ one
has
\begin{eqnarray}
\int\limits_{h_1(x)}^x \omega\, >\, N,\label{int}
\end{eqnarray}
where the integral is calculated along the path $t\mapsto h_{1-t}(x)\in X$,
$t\in [0,1]$.
\end{definition}

Intuitively, an $N$-movable subset can be continuously deformed inside $X$ such
that each of its points is winding around yielding a large quantity (\ref{int})
measured in terms of the form $\omega$.

It is easy to see that (assuming that $\xi\not=0$ and the space $X$ is
connected) any subset $A\subset X$ such that the inclusion $A\to X$ is
null-homotopic is $N$-movable for any $N$ with respect to any closed 1-form
$\omega$ representing the class $\xi$.

\begin{example}\label{ex1}
{\rm Assume that $X=Y\times S^1$ and $\xi\in H^1(X;\R)$ is such that
$\xi|_{S^1}\not=0$. Then for any integer $N$ the total space $A=X$ is
$N$-movable with respect to any closed 1-form $\omega$ representing class
$\xi$. A homotopy $h_t: X\to X$ as above can be described as follows. Identify
$S^1\subset \C$ with the set of complex numbers having norm one. Then for $y\in
Y$ and $z\in S^1$ one sets
$$h_t(y,z) = (y, e^{iMt}\cdot z), \quad t\in [0,1]$$
with suitable real $M$.}
\end{example}
\begin{definition}\label{deff3}
Fix a closed 1-form $\omega$ representing $\xi$. The number $\cat(X,\xi)$ is
the minimal integer $k$ with the property that  for any $N>0$ there exists an
open cover $F, F_1, \dots, F_k\subset X$ such that each inclusion $F_j\to X$ is
null-homotopic and such that $F$ is $N$-movable with respect to $\omega$.
\end{definition}

\begin{example}\label{ex2}
{\rm Assume that $X=Y\times S^1$ and $\xi\in H^1(X;\R)$ is such that
$\xi|_{S^1}\not=0$. Then $\cat(X,\xi)=0$, compare Example \ref{ex1}.}
\end{example}

It is known that $\cat(X,\xi)$ is homotopy invariant, see Lemma 10.12 in
\cite{farbe3}. In particular, it is independent of the choice of $\omega$, see
\cite{farbe3}, page 166.

The next result gives a lower bound for $\cat(X,\xi)$ in terms of cohomological
information.

Denote by $p:\tilde X\to X$ the covering corresponding to $\ker \xi$.

\begin{theorem}\label{thm1}
Suppose that $ L_0, L_1, \dots, L_k\in \mathcal V_\xi$ are complex flat vector
bundles with the following properties: (1) $L_0$ does not belong to the support
$\Supp(X,\xi)$; (2) For a homology class
$$z\in H_d(\tilde X;\C)$$ and some cohomology classes
$$v_i\in H^{d_i}(X;
L_i), \quad i=0, 1, \dots, k,$$ where
\begin{eqnarray}\label{pos}d_i>0,\quad \mbox{for
$i=1, \dots, k$}\end{eqnarray}
the cup-product $ v_0\cup v_1\cup \dots \cup
v_k\in H^d(X;L)$ evaluates nontrivially on the homology class
$$p_\ast(z)\in H_d(X;L^\ast)$$ i.e.
\begin{eqnarray}\label{nonzero}
\langle v_0\cup \dots \cup v_k, p_\ast(z)\rangle \not=0.
\end{eqnarray} Here
$$z\in H_d(\tilde X;\C), \quad d= \sum_{i=0}^k d_i, \quad L= L_0\otimes L_1\otimes \dots \otimes
L_k$$ and $L^\ast$ is the dual of $L$. Then
\begin{eqnarray}\label{second}
\cat(X,\xi) \, > \,  k.
\end{eqnarray}
\end{theorem}

\begin{remark} {\rm Theorem \ref{thm1} is meaningful for $k=0$ as well; in this case it
gives $\cat(X,\xi)>0$. }
\end{remark}

\begin{remark} {\rm  In the case $\xi=0$ (i.e. when one deals with functions)
the statement above turns into the usual cup-length
estimate for the Lusternik - Schnirelmann category $\cat(X)$, see \S \ref{sec5}
below.}
\end{remark}

\begin{proof}[Proof of Theorem \ref{thm1}]
Assume that (\ref{second}) is false, i.e. $\cat(X,\xi)\leq k$. Let $\omega$ be
a continuous closed 1-form on $X$ representing $\xi$.

Then for any $N>0$ there exists an open cover $F, F_1, \dots, F_k\subset X$
such that each inclusion $F_j\to X$ is null-homotopic and such that $F$ is
$N$-movable with respect to $\omega$.

Fix a singular cycle $c$ in $\tilde X$ representing the class $z\in H_d(\tilde
X;\C)$. One may find a compact polyhedron $K\subset \tilde X$ containing $c$.

Let us show that there exists a neighborhood of infinity $U\subset \tilde X$
with the following property: any homology class in $H_\ast(K;\C)$, which is
homologous in $\tilde X$ to a cycle lying in $U$, is movable to infinity of
$\tilde X$ with respect to $\xi$. Indeed, start with an arbitrary neighborhood
of infinity $N$ satisfying Lemma 3 of \cite{FS}. For a covering translation
$g:\tilde X \to \tilde X$ the intersection
$$V_g=\im[H_\ast(gN)\to H_\ast(\tilde X)]\cap \im[H_\ast(K)\to H_\ast(\tilde X)]$$
(we use homology with $\C$ coefficients) is a finite dimensional complex vector
space. Therefore one may find a covering translation $g_0:\tilde X\to \tilde X$
such that for any covering translation $g$ one has $V_{g_0}\subset V_g$. Then
clearly the neighborhood of infinity $U=g_0N$ satisfies the above requirement.

One has $p^\ast\omega =df$ where $f: \tilde X\to \R$ is a continuous function.
Then $f(K)\subset [a, b]$ and $U\supset f^{-1}(-\infty, c)$ where $c<a<b$.

Pick a number $N> b-c$ and apply Definition \ref{deff3}. We obtain an open
cover $F\cup F_1\cup \dots\cup F_k=X$ where $F$ is $N$-movable with respect to
$\omega$ and each inclusion $F_j\to X$ is null-homotopic, $j=1, 2, \dots, k$.

Find subsets $A'\subset A\subset F$ and $B'\subset B\subset F_1\cup F_2\cup
\dots\cup F_k$ with the following properties: (i) $A'$ and $B'$ are open and
cover $X$, i.e. $A'\cup B'=X$; (ii) $A$ and $B$ are compact sub-polyhedra of
$X$.

The restriction of the cup-product $v = v_1\cup \dots \cup v_k$ on the set
$F_1\cup \dots \cup F_k$ vanishes\footnote{Here we use  (\ref{pos}) and the
assumption that $F_j\to X$ is null-homotopic, $j=1, \dots k$.} and so $v$ can
be realized by a singular cochain $g$ vanishing on all singular simplices lying
entirely in $B'$. The cochain $g$ takes values in the local system $L_1\otimes
\dots\otimes L_k$.

We may assume that the topology of $X$ is given by a metric $d$. Let $\epsilon
>0$ be the Lebesgue number of the open cover $A'\cup B'$.

Subdivide the singular chain $c$ representing $z$ such that it is a linear
combination of finitely many singular simplices, each of diameter $<\epsilon$.

From (\ref{nonzero}) we have \begin{eqnarray}\label{nonzero1} \langle v_0,
p_\ast(p^\ast(v)\cap z)\rangle \not= 0.\end{eqnarray} The class
$$z_0=p^\ast(v)\cap z\in H_\ast(\tilde X;\C)$$ is represented by the singular
cycle $p^\ast(g)\cap c$ (see \cite{Sp}, chapter 5, \S 6) having support in $K$.
From the construction above, we see that the projected cycle
$p_\ast(p^\ast(g)\cap c)$ lies entirely in $A'$. We obtain (because of our
choice of $N$ and $U$) that the homology class $z_0$ is movable to infinity of
$\tilde X$ with respect to $\xi$. But this contradicts Theorem \ref{cor4} since
$\langle v_0, p_\ast(z_0)\rangle\not=0$ and the bundle $L_0$ is assumed to be
generic, i.e. $L_0\notin \Supp(X,\xi)$.
\end{proof}

\section{Transcendental line bundles}\label{sectrans}

In this section we improve Theorem \ref{thm1} by showing that under certain
conditions one may avoid mentioning explicitly the homology class $z\in
H_\ast(\tilde X)$ in the statement.

First we recall our notations.
Let $\xi\in H^1(X;\R)$ be a real cohomology
where $X$ is a finite cell complex. $\xi$ determines a homomorphism
$H_1(X;\Z)\to \R$; we denote by $\ker(\xi)$ its kernel. We set
$$H = H_1(X;\Z)/\ker ( \xi).$$
It is a free abelian group of finite rank $r=\rk \xi$. Any flat line bundle
$L\in \V_\xi$ determines a monodromy homomorphism (\ref{monodromy}).

\begin{definition}\label{transcendental}
We say that a bundle $L\in \V_\xi$  is algebraic if the monodromy homomorphism
${\rm Mon}_L: \Z[H]\to \C$ has nontrivial kernel. We say that $L$ is
transcendental\footnote{Perhaps it would be better to call such bundles {\it
$\xi$-algebraic} and {\it $\xi$-transcendental} as both properties depend on
the class $\xi$ (they actually depend only on $\ker \xi$).} if ${\rm Mon}_L:
\Z[H]\to \C$ is injective.
\end{definition}

If $t_1, \dots, t_r\in H$ is a basis and if $a_i\in \C$ denotes the monodromy
of $L$ along $t_i$, i.e. $a_i={\rm Mon}_L(t_i)$, then $L$ is algebraic iff
there exist a nontrivial Laurent polynomial equation with integral coefficients
$q(t_1, \dots, t_r)$ such that $q(a_1, \dots, a_r)=0$.

There exist countably many nonzero Laurent polynomials $q$ with integral
coefficients and for each such $q$ the set of solutions $q(a_1, \dots, a_r)=0$
is nowhere dense in $\V_\xi$. Since $\V_\xi=(\C^\ast)^r$ is homeomorphic to a
complete metric space we obtain:

\begin{lemma}\label{bair} The set of all transcendental bundles $L\in \V_\xi$ is of Baire category 2. In
particular, the set of transcendental $L\in \V_\xi$ is dense in the variety
$\V_\xi$.
\end{lemma}

\begin{lemma}\label{lomeshane} The dimension of the vector space $H^q(X;L)$ is constant on the
set of transcendental flat line bundles $L\in \V_\xi$. In other words, for any
two transcendental $L, L'\in \V_\xi$ one has $$\dim H^q(X;L)=\dim H^q(X;L').$$
\end{lemma}
\begin{proof} The monodromy homomorphism ${\rm {Mon}}_L: \Z[H]\to \C$ defines a
left $\Z[H]$-module structure $\C_L$ on $\C$ and by definition
\begin{eqnarray}\label{deftwisted}
H^q(X;L)=H^q(\Hom_{\Z[H]}(C_\ast(\tilde X), \C_L)) \end{eqnarray}
 where
$C_\ast(\tilde X)$ is the cellular chain complex of the covering $\tilde X\to
X$ corresponding to $\ker(\xi)$. If $L$ is transcendental then ${\rm {Mon}}_L$
gives a field extension $Q(H)\to \C$ where $Q(H)$ is the field of fractions of
$\Z[H]$. We obtain therefore (using finiteness of $C_\ast(\tilde X)$ over
$\Z[H]$ and (\ref{deftwisted})):
\begin{eqnarray*}
H^q(X;L) &\simeq & H^q(\Hom_{\Z[H]}(C_\ast(\tilde X);Q(H)))\otimes_{Q(H)}\C_L
\\ &\simeq & H^q(X;Q(H))\otimes_{Q(H)}\C_L. \end{eqnarray*}
This implies that for any transcendental $L$ one has
\begin{eqnarray}\label{nov}
\dim_\C H^q(X;L) \, =\, \dim_{Q(H)}H^q(X;Q(H))
\end{eqnarray}
and the right hand side\footnote{The common value $\dim_\C H^q(X;L)$ for $L$
transcendental which appears in (\ref{nov}) equals the Novikov -- Betti number
$b_q(\xi)$, see \cite{noviko} and \cite{farbe3}, Proposition 1.30.} is
independent of $L$.
\end{proof}

\begin{proposition}\label{notin} Let $X$ be a finite complex and $\xi\in H^1(X;\R)$. If the flat
line bundle $L\in
\V_\xi$ is transcendental then $L\notin \Supp(X,\xi)$.
\end{proposition}

\begin{proof} If $L$ is transcendental then the monodromy homomorphism ${\rm
{Mon}_L}$ can be decomposed into $\Z[H]\to Q(H)\to \C$ (where we use the
notations introduced in the proof of the previous lemma) and hence the
homomorphism $p_\ast: H_q(\tilde X;\Z) \to H_q(X;L)$ can be decomposed into
$H_q(\tilde X;\Z)\to Q(H)\otimes_{\Z[H]}H_q(\tilde X;\Z)\to H_q(X;L)$. This
shows that all $\Z[H]$-torsion classes $z\in H_q(\tilde X;\Z)$ satisfy
$p_\ast(z)=0$.

Let $T_q\subset H_q(\tilde X;\Z)$ be the subgroup of all $\Z[H]$-torsion
classes. We claim that $\C[H]$-torsion of $H_q(\tilde X;\C)$ coincides with
$\C[H]\otimes_{\Z[H]}T_q$. Together with the remark of the previous paragraph
this would imply that for any $\C[H]$-torsion class $z\in H_q(\tilde X;\C)$ one
has $p_\ast(z)=0$, i.e. $L\notin \Supp(X,\xi)$.

We have the exact sequence
\begin{eqnarray}\label{exact}
\qquad 0\to \C[H]\otimes T_q \to \C[H]\otimes H_q(\tilde X;\Z) \to \C[H]\otimes
M\to 0
\end{eqnarray}
where $M$ is $H_q(\tilde X;\Z)/T_q$ and the tensor product is over $\Z[H]$. The
sequence (\ref{exact}) is exact since $\C[H]$ is flat as a $\Z[H]$-module.
Using \cite{B}, exercise 20 from chapter 1, \S 2 on page 46, we find that
$\C[H]\otimes_{\Z[H]}M$ has no $\C[H]$-torsion. The sequence (\ref{exact})
implies now that $\C[H]\otimes_{\Z[H]}T_q$ coincides with the $\C[H]$-torsion
of $H_q(\tilde X;\C)=\C[H]\otimes_{\Z[H]}H_q(\tilde X;\Z)$.
\end{proof}

The authors are thankful to Holger Brenner for valuable comments leading to
Proposition \ref{notin}.

\begin{proposition}\label{propnew}
Assume that $L\in \V_\xi$ is transcendental. Let $v\in H^q(X;L)$ be a nonzero
cohomology class. Then there exists a homology class $z\in H_q(\tilde X;\Z)$
such that $\langle v, p_\ast(z)\rangle \not=0$. Here
\begin{eqnarray}\label{proj1}
p_\ast: H_q(\tilde X;\Z)\to H_q(X;L^\ast)
\end{eqnarray}
denotes homomorphism similar to (\ref{proj}). Moreover, the kernel of
(\ref{proj1}) coincides with the set of $\Z[H]$-torsion classes in $H_q(\tilde
X;\Z)$.
\end{proposition}
\begin{proof}
Obviously, if the cohomology class $v$ is nonzero then there exists a homology
class $z'\in H_q(X;L^\ast)$ with $\langle v, z'\rangle \not=0$. Here $L^\ast\in
\V_\xi$ denotes the dual bundle to $L$.

If $L$ is transcendental then the dual bundle $L^\ast$ is transcendental as
well, i.e. ${\rm Mon}_{L^\ast}: \Z[H]\to \C$ is injective. Let $\C_{L^\ast}$
denote the field of complex numbers $\C$ viewed as a right $\Z[H]$-module via
the ring homomorphism ${\rm Mon}_{L^\ast}$. We want to show that
\begin{eqnarray}\label{tensor1}
H_q(X;L^\ast) \simeq \C_{L^\ast}\otimes_{\Z[H]}H_q(\tilde X;\Z).
\end{eqnarray}
Assuming (\ref{tensor1}) we would be able to argue that the class $z'$ can be
represented as a finite sum $z'=\sum c_i p_\ast(z_i)$ where $c_i\in \C$ and
$z_i\in H_q(\tilde X;\Z)$. Since $\langle v, z'\rangle \not=0$ the number
$\langle v, p_\ast(z_i)\rangle$ is nonzero for some $i$ which implies our
statement.

Observe that $\Z[H]$ is a Laurent polynomial ring. Denote by $Q(H)$ its field
of fractions. The ring homomorphism ${\rm Mon}_{L^\ast}:\Z[H]\to \C_{L^\ast}$
extends to a field embedding $Q(H)\to \C_{L^\ast}$. Let $C_\ast(\tilde X)$ be
the cellular chain complex of $\tilde X$ with integral coefficients. The
homomorphism (\ref{proj1}) is induced by the chain map
\begin{eqnarray*}
C_\ast(\tilde X) \to \C_{L^\ast}\otimes_{\Z[H]}C_\ast(\tilde X).
\end{eqnarray*}
The latter can be decomposed as
\begin{eqnarray}\label{rational}
\begin{array}{ll}
C_\ast(\tilde X) \to Q(H)\otimes_{\Z[H]}C_\ast(\tilde X) \to
\C_{L^\ast}\otimes_{\Z[H]}C_\ast(\tilde X) =\\ \\ = \C_{L^\ast}\otimes
_{Q(H)}(Q(H)\otimes_{\Z[H]}C_\ast(\tilde X)).
\end{array}
\end{eqnarray}
The left map in (\ref{rational}) is a localization and hence it induces
localization on homology
\begin{eqnarray*}
H_q(\tilde X;\Z)\to Q(H)\otimes_{\Z[H]}H_q(\tilde X;\Z).
\end{eqnarray*}
The right map in (\ref{rational}) is induced by a field extension $Q(H)\to
\C_{L^*}$. Note that $\C_{L^\ast}$ viewed as a $\Z[H]$-module is a direct sum
of infinite number of copies of $Q(H)$. Hence the right homomorphism in
(\ref{rational}) induces
\begin{eqnarray*}
Q(H)\otimes_{\Z[H]}H_q(\tilde X;\Z) \to \C_{L^\ast}\otimes_{\Z[H]}H_q(\tilde
X;\Z).
\end{eqnarray*}
Therefore (\ref{proj1}) coincides with $H_q(\tilde X;\Z)\to
\C_{L^\ast}\otimes_{\Z[H]}H_q(\tilde X;\Z)$ proving (\ref{tensor1}).
\end{proof}

Next we state an improved version of Theorem \ref{thm1} taking into account
Proposition \ref{propnew}.

\begin{theorem}\label{thm11}
Suppose that $ L_0, L_1, \dots, L_k\in \mathcal V_\xi$ are complex flat vector
bundles with the following properties: (1) $L_0$ does not belong to the support
$\Supp(X,\xi)$; (2) For some cohomology classes $v_i\in H^{d_i}(X; L_i)$, $i=0,
1, \dots, k,$ where $d_i>0$ for $i=1, \dots, k$ the cup-product
\begin{eqnarray} v_0\cup v_1\cup
\dots \cup v_k\not=0\in H^d(X;L)\end{eqnarray}
 is nontrivial; here
$d= \sum_{i=0}^k d_i$ and $ L= L_0\otimes L_1\otimes \dots \otimes L_k$. (3)
$L$ is transcendental, i.e. the monodromy homomorphism ${\rm Mon}_L: \Z[H]\to
\C$ is injective. Then
\begin{eqnarray}\label{second2}
\cat(X,\xi) \, > \,  k.
\end{eqnarray}
\end{theorem}
\begin{proof}
It follows by combining Theorem \ref{thm1} with Proposition \ref{propnew}.
\end{proof}

As a useful special case of Theorem \ref{thm11} and Proposition \ref{notin} we
mention the following statement:

\begin{theorem}\label{perfect}
Let $X$ be a finite cell complex and $\xi\in H^1(X;\R)$. Let $L\in \V_\xi$ be
transcendental. Assume
 that there exist
cohomology classes $v_0\in H^{d_0}(X;L)$ and $v_i\in H^{d_i}(X;\C)$ where $i=1,
\dots, k$ such that $d_i>0$ for $i\in \{1, \dots, k\}$ and the cup-product
\begin{eqnarray}
v_0\cup v_1\cup\dots\cup v_k\, \not=0 \, \in H^\ast(X;L)
\end{eqnarray}
is nontrivial. Then $\cat(X,\xi)>k$.
\end{theorem}

Theorem \ref{perfect}, combining simplicity with remarkable efficiency, has a
very satisfying statement. We view this theorem as being the main result of the
paper. In the following sections we test this theorem in many specific
examples. Besides, we compare theorems of this section with cohomological lower
bounds for $\cat(X,\xi)$ obtained earlier in \cite{farbe3}, Th. 6.1 and in
\cite{farbe2}, Th. 10.23.

\section{The notion of cup-length $\cl(X,\xi)$}

In view of Theorem \ref{perfect} we introduce the following notation.

Let $X$ be a finite cell complex and $\xi\in H^1(X;\R)$. We denote by
$\cl(X,\xi)$ the maximal integer $k\geq 0$ such that Theorem \ref{perfect}
could be applied to $(X,\xi)$; if Theorem \ref{perfect} is not applicable (i.e.
if $H^\ast(X;L)=0$ for any transcendental $L\in \V_\xi$, compare Lemma
\ref{lomeshane}) we set
$$\cl(X,\xi)=-1.$$
Hence, \begin{eqnarray}\cl(X,\xi)\, \in \, \{-1, 0, 1, \dots\}.
\end{eqnarray}

In other words, $\cl(X,\xi)\geq k$ where $k\geq 0$ iff there exists a
transcendental flat line bundle $L\in \V_\xi$ and there exist cohomology
classes $v_0\in H^{d_0}(X;L)$ and $v_i\in H^{d_i}(X;\C)$ where $i=1, \dots, k$
and $d_i>0$ for $i\in \{1, \dots, k\}$ such that the cup-product
\begin{eqnarray}\label{cupprod}
v_0\cup v_1\cup\dots\cup v_k\, \not=0 \, \in H^\ast(X;L)
\end{eqnarray}
is nontrivial.

Note that for $\xi=0$ the number $\cl(X,\xi)$ coincides with the usual
cup-length $\cl(X)$; recall that the later is defined as the largest integer
$r$ such that there exist cohomology classes $u_i\in H^{d_i}(X;\C)$ where $i=1,
\dots, r$ of positive degree such that their cup-product $ u_1\cup\dots\cup
u_k\, \not=0 \, \in H^\ast(X;L) $ is nontrivial. Indeed, in the case $\xi=0$
the trivial bundle $L=\C$ is not in $\Supp(X,\xi)$ and therefore one may take
$L=\C$ and $v_0=1\in H^0(X;\C)$.

One can restate Theorem \ref{perfect} as follows:

\begin{theorem}\label{rephrase} One has $\cat(X,\xi)\geq \cl(X,\xi)+1$.
\end{theorem}

The next useful Lemma suggests several different ways to characterize the
number $\cl(X,\xi)$. This Lemma plays an important role in the sequel.

\begin{lemma}\label{equivalent} Let $X$ be a finite cell complex and $\xi\in H^1(X;\R)$.
The following statements regarding an integer $k\geq 0$ are equivalent:
\begin{enumerate}
\item[{\rm (A)}] $\cl(X,\xi)\geq k$;

\item[{\rm (B)}] There exists cohomology class $v_0\in H^{d_0}(X;L)$ where
$L\in \V_\xi$ is a transcendental flat line bundle and there exist $k$ integral
cohomology classes $v_i\in H^{d_i}(X;\Z)$ where $i=1, \dots, k$ and $d_i>0$ for
$i\in \{1, \dots, k\}$, such that the cup-product (\ref{cupprod}) is
nontrivial.

\item[{\rm (C)}] Let $H$ denote $H_1(X;\Z)/\ker \xi$; it is a free abelian
group of finite rank. Let $Q(H)$ denote the field of fractions of the group
ring $\Z[H]$. Then there exist cohomology classes $w_0\in H^{d_0}(X;Q(H))$ (the
latter denotes homology with twisted coefficients) and $v_i\in H^{d_i}(X;\Z)$
where $d_i>0$ for $i=1, \dots, k$ such that the cup-product
\begin{eqnarray*}
w_0\cup v_1\cup \dots \cup v_k\not=0\in H^\ast(X;Q(H))
\end{eqnarray*}
is nontrivial.

\item[{\rm (D)}] For any transcendental flat line bundle $L\in \V_\xi$ there
exists cohomology class $v_0\in H^{d_0}(X;L)$ such that the cup-product
(\ref{cupprod}) is nontrivial for some integral cohomology classes $v_i\in
H^{d_i}(X;\Z)$ with $d_i>0$ for $i\in \{1, \dots, k\}$.

\item[{\rm (E)}] For any transcendental flat line bundle $L\in \V_\xi$ there
exist cohomology classes $v_0\in H^{d_0}(X;L)$ and $v_i\in H^{d_i}(X;\C)$ where
$i=1, \dots, k$ and $d_i>0$ for $i\in \{1, \dots, k\}$ such that the
cup-product (\ref{cupprod}) is nontrivial.
\end{enumerate}
\end{lemma}
\begin{proof}
Let us show that (A)$\implies$(B). Fix a transcendental bundle $L\in \V_\xi$,
and $v_0\in H^\ast(X;L)$ such that (\ref{cupprod}) is nontrivial for some
$v_i\in H^{d_i}(X;\C)$ with $d_i>0$. Consider now the cup-products
\begin{eqnarray}\label{function}
v_0\cup v'_1\cup\dots \cup v'_k
\end{eqnarray}
with arbitrary integral cohomology classes $v'_i\in H^{d_i}(X;\Z)$; here the
degrees $d_i$ are assumed to be fixed. (\ref{function}) is a multi-linear
function of the classes $v'_i$. Since the integral classes generate
$H^{d_i}(X;\C)$ over $\C$ we obtain that (\ref{function}) must be nonzero for
some choice of classes $v'_i$, i.e. (B) holds.

Now we show that (B)$\implies$(C). Fix $L\in \V_\xi$ and the classes $v_0\in
H^{d_0}(X;L)$ and $v_i\in H^{d_i}(X;\Z)$ satisfying conditions described in
(B). The monodromy homomorphism ${\rm {Mon}}_L: \Z[H]\to \C$ is an injective
ring homomorphism, it extends to the field of fractions $Q(H)\to \C$. The image
of the induced homomorphism on cohomology
\begin{eqnarray}\label{flat}
\psi: H^{d_0}(X;Q(H)) \to H^{d_0}(X;L)
\end{eqnarray}
generates $H^{d_0}(X;L)$ over $\C$ and (\ref{flat}) is injective (for reasons
mentioned in the proof of Lemma \ref{lomeshane}). Fix cohomology classes
$v_i\in H^\ast(X;\Z)$ where $i=1, \dots, k$. For a cohomology class $w_0\in
H^{d_0}(X;Q(H))$ the function
\begin{eqnarray}\label{com} \psi(w_0)\cup v_1\cup \dots\cup v_k= \psi(w_0\cup
v_1\cup \dots\cup v_k)\, \in \, H^\ast(X;L)
\end{eqnarray} extends to a $\C$-linear function of $v_0\in H^{d_0}(X;L)$
$$v_0\mapsto v_0\cup v_1\cup \dots\cup v_k \in H^\ast(X;L).$$
If this function is nonzero then it may not vanish on the image of $\psi$, i.e.
(C) holds.

The implication (C)$\implies$(D) follows from injectivity of homomorphism
(\ref{flat}) and from (\ref{com}).

Implications (D)$\implies$(E) and (E)$\implies$(A) are obvious. This completes
the proof.
\end{proof}

\begin{lemma}\label{clgeq}
Assume that $X$ and $Y$ are path connected finite cell complexes and $\xi\in
H^1(X\times Y;\R)$.
 Then
\begin{eqnarray}\label{cl}
\cl(X\times Y, \xi)\geq \cl(X,\xi|_X) + \cl(Y, \xi|_Y).
\end{eqnarray}
\end{lemma}
\begin{proof} Denote $\cl(X, \xi|_X)=k$ and $\cl(Y, \xi|_Y)=r$.
Any flat line bundle $L$ over $X\times Y$ has the form $L_1\boxtimes L_2$
(exterior tensor product) where $L_1$ and $L_2$ are flat line bundles over $X$
and $Y$ respectively. Note that if $L$ lies in the variety
$\V_\xi=\Hom(H_1(X\times Y;\Z)/\ker (\xi),\C^\ast)$ then $L_1$ and $L_2$ are
obtained by restrictions and hence $L_1\in \V_{\xi|_X}$ and $L_2\in
\V_{\xi|_Y}$. We will use equivalence between (A), (E) and (F) of Lemma
\ref{equivalent}. Fix a transcendental bundle $L=L_1\boxtimes L_2\in \V_\xi$
over $X\times Y$. Then both $L_1$ and $L_2$ are transcendental. Find classes
$v_0\in H^\ast(X;L_1)$, $v_1, \dots, v_k\in H^\ast(X;\C)$, $u_0\in
H^\ast(Y;L_2)$, $u_1, \dots, u_r\in H^\ast(Y;\C)$ such that $v_0\cup v_1\cup
\dots\cup v_k\not=0$ and $u_0\cup u_1\cup \dots\cup u_r\not=0$. Now we have
cohomology classes $v_0\times u_0\in H^\ast(X\times Y;L)$ and $v_i\times 1,
1\times u_j\in H^\ast(X\times Y;\C)$ and the product
\begin{eqnarray*}
(v_0\times u_0) \, \cup \, \prod_{i=1}^k (v_i\times 1) \, \cup \, \prod_{j=1}^r
(1\times u_j)\not=\, \, 0\, \in\,  H^\ast(X\times Y;L)
\end{eqnarray*}
is nonzero. Here we use the K\"unneth formula which states
\begin{eqnarray}
H^\ast(X\times Y; L_1\boxtimes L_2) \simeq H^\ast(X;L_1) \otimes H^\ast(Y;L_2).
\end{eqnarray}
This proves (\ref{cl}).
\end{proof}

\section{First examples}\label{sec5}

In this section we test Theorem \ref{perfect} in a number of simple examples.

\subsection{}
First consider the case $\xi=0$. We know that for $\xi=0$ the number
$\cat(X,\xi)$ coincides with the classical LS category $\cat(X)$, see
\cite{farbe3}, Example 10.8. Let us examine what Theorem \ref{perfect} gives in
this case. The variety $\V_\xi$ has only one point -- the trivial flat line
bundle $\C$ over $X$. The support $\Supp(X,\xi)=\emptyset$ is always empty for
$\xi=0$. We may therefore take $v_0=1\in H^0(X;\C)$ applying Theorem
\ref{perfect}. Thus, we see that Theorem \ref{perfect} claims in the special
case $\xi=0$ that if there exist cohomology classes $v_1, \dots, v_k\in
H^{>0}(X;\C)$ with $v_1\cup\dots\cup v_k\not=0$ then $\cat(X)>k$. This claim is
the classical cup-length estimate for the LS category $\cat(X)$.

\subsection{}
Note that if $\xi\not=0$ and $X$ is connected then $H^0(X;L)=0$ for any
nontrivial $L\in \V_\xi$. Note also that the trivial flat line bundle $\C\in
\V_r$ always lies in the support $\Supp(X,\xi)$ for $\xi\not=0$, see Remark
\ref{rmk1} and it is never transcendental. Therefore the degree of the class
$v_0$ (which appears in Theorem \ref{thm1}) in the case $\xi\not=0$ must be
positive. Hence for $\xi\not=0$ the number $k$ in Theorem \ref{thm11} satisfies
$k \leq \dim X-1$. This explains why Theorem \ref{perfect} cannot give
$\cat(X,\xi) \geq \dim X+1$ for $\xi\not=0$.

Inequality (10.8) in \cite{farbe3} yields
\begin{eqnarray}\label{lesss}
\cat(X,\xi)\leq \cat(X)-1 \leq \dim X
\end{eqnarray}
assuming that $X$ is connected and $\xi\not=0$. This is consistent with the
remark of the previous paragraph.

\subsection{}
The following example shows that (\ref{lesss}) can be satisfied as an equality
i.e. that $\cat(X,\xi)=\dim X$ is possible\footnote{Note that if $X$ is a
closed smooth manifold then $\cat(X,\xi)< \dim X$, see \S \ref{sec6}.}.
Consider the bouquet $X=Y\vee S^1$ where $Y$ is a finite polyhedron, and assume
that the class $\xi\in H^1(X;\R)$ satisfies $\xi|_Y=0$ and $\xi|_{S^1}\not=0.$
We know that in this case
\begin{eqnarray}\label{minus} \cat(X,\xi) =
\cat(Y)-1,\end{eqnarray} see Example 10.11 from \cite{farbe3}.

We are going to apply Theorem \ref{perfect}. The variety $\V_\xi$ in this case
coincides with the set $\C^\ast=\C-\{0\}.$ The support $\Supp(X,\xi)$ contains
in this case only the trivial line bundle. $L\in \V_\xi$ is transcendental if
the monodromy along the circle $S^1$ is a transcendental complex number. For
any $L\in \V_\xi$ the restriction $L|_Y$ is trivial and the restriction
homomorphism $H^i(X;L)\to H^i(Y;\C)$ is onto.

Suppose that the cohomological cup-length of $Y$ with $\C$-coefficients equals
$\ell$, i.e. there exist cohomology classes of positive degree $u_0, u_1,
\dots, u_{\ell-1}\in H^{>0}(Y;\C)$ such that the product $u_0\cup \dots\cup
u_{\ell-1}\not=0$ is nonzero. By the above remark, for a nontrivial $L\in
\V_\xi$ we obtain cohomology classes $v_0\in H^\ast(X;L)$ and $v_1, \dots,
v_{\ell-1}\in H^{>0}(X;\C)$ such that $v_i|_Y=u_i$. Hence $v_0\cup v_1\cup
\dots\cup v_{\ell-1} \not=0\in H^\ast(X; L)$. By Theorem \ref{perfect} we
obtain $\cat(X,\xi) > \ell-1$ which is equivalent (taking into account
(\ref{minus})) to $\cat(Y)>\ell$. The last inequality is the classical
cup-length estimate for the usual category.

\subsection{}
Comparing the previous example with Example 6.6 in \cite{farbe2} we find that
Theorem \ref{thm1} and Theorem \ref{perfect} are stronger than the
cohomological lower bounds given in \cite{farbe2}, even in the case of rank one
classes $\xi$. The reason for this is that in Theorem \ref{thm1} we have only
one generic bundle (compared with two in Theorem 6.1 from \cite{farbe2}) but
non-vanishing of the product $v_0\cup \dots\cup v_k$ is understood in a
stronger sense.

\subsection{} \label{specific}
Let us now consider a very specific example: $X=T^2\vee S^1$. In this case
$H^1(X;\R)=\R^3$ and we describe $\cat(X,\xi)$ as function of $\xi\in
H^1(X;\R)=\R^3$. We denote by $\ell\subset \R^3$ the set of all classes $\xi$
such that $\xi|_{T^2}=0$. Clearly $\ell$ is a line through the origin in
$\R^3$. We claim that:
\begin{eqnarray}\label{cases}
\cat(X,\xi)=\left\{ \begin{array}{ll} 1, & \mbox{if $\xi\in \R^3 -\ell$},\\
2, & \mbox{if $\xi\in \ell-\{0\}$},\\
3, & \mbox{if $\xi=0$}.
\end{array}
\right.
\end{eqnarray}
Indeed, consider first the case $\xi\notin \ell$, i.e. $\xi|_{T^2}\not=0$. Let
us show that $\cat(X,\xi)\leq 1$. Denote $p=T^2\cap S^1$ and let $q\in S^1$ be
a point distinct from $p$. Set $F=X-\{q\}$ and $F_1=S^1-\{p\}$. Then $F\cup
F_1=X$ is an open cover of $X$ with $F_1\to X$ null-homotopic and with $F$
being $N$-movable in $X$ for any $N$ (assuming that $\xi|_{T^2}\not=0$; this
follows from homotopy invariance of $\cat(X,\xi)$ and from Examples \ref{ex1}
and \ref{ex2}.

Since $\cat(X,\xi)=0$ would imply $\chi(X)=0$ by Theorem \ref{ep} stated below
we obtain that $\cat(X,\xi)>0$ for any $\xi$ (as $\chi(X)=-1\not=0$). This
proves the first line of (\ref{cases}).

If $\xi\in \ell-\{0\}$ we apply the result of Example 3.5 from \cite{farbe2}
which gives $\cat(X,\xi)=\cat(T^2)-1=2$.

For $\xi=0$ we easily find $\cat(X,\xi)=\cat(X)=3$.

\section{A controlled version of $\cat(X,\xi)$}\label{control}

In this section we introduce a new controlled version of $\cat(X,\xi)$ which
has some advantages, for example it behaves better under products.

Let $\omega$ be a continuous closed 1-form on a finite cell complex $X$. Let
$\xi=[\omega]\in H^1(X;\R)$ be the cohomology class represented by $\omega$. We
refer to \cite{farbe2} for the formalism of closed 1-forms on topological
spaces.

\begin{definition}\label{cmovable}
Let $N$ and $C$ be two real numbers. A subset $A\subset X$ is $N$-movable with
respect to $\omega$ with control $C$ if there exists a continuous homotopy
$h_t:A\to X,$ $t\in [0,1],$ such that (1) $h_0:A\to X$ is the inclusion; (2)
for any point $x\in A$ one has
\begin{eqnarray}
\int\limits_{x}^{h_1(x)} \omega\,<\,- N,\label{cint}
\end{eqnarray}
where the integral is calculated along the path $t\mapsto h_{t}(x)\in X$, $t\in
[0,1]$ and (3) for any point $x\in A$ and for any $t\in [0,1]$ one has
\begin{eqnarray}
\int\limits_{x}^{h_t(x)} \omega\,\leq\,C.\label{ccint}
\end{eqnarray}
\end{definition}

Geometrically, an $N$-movable subset with control $C$ can be continuously
deformed inside $X$ such that each of its points is winding around yielding a
small quantity (\ref{cint}) measured in terms of the form $\omega$ and such
that for all times $t$ the integral (\ref{ccint}) remains controlled, i.e. it
is smaller than a fixed quantity $C$.

\begin{example} {\rm Let $X$ be a closed smooth manifold and let $\omega$ be a
smooth closed 1-form on $X$ having no zeroes. Find a vector field $v$ on $X$
such that $\omega(v)<0$. Then the whole space $X$ is $N$-movable with control
$C=0$ for any $N>0$. The homotopy $h_t: X\to X$ as in Definition \ref{cmovable}
can be easily constructed using the flow generated by $v$.}
\end{example}

\begin{example}\label{ex12}
{\rm Assume that $X$ is a path connected polyhedron and the cohomology class of
$\omega$ is nonzero $\xi=[\omega]\not=0\in H^1(X;\R)$. For any subset $F\subset
X$ such that the inclusion $F\to X$ is homotopic to a constant map $F\to X$
there exists a constant $C>0$ so that for any $N>0$ the set $F$ is $N$-movable
with respect to $\omega$ with control $C$. Indeed, find a closed loop $\gamma$
in $X$ such that $\langle \xi, [\gamma]\rangle =\int_\gamma \omega <0.$ Define
the homotopy $h_t:F\to X$ as concatenation of contraction of $F$ to a point in
$X$, then moving this point towards the loop $\gamma$ and finally traversing
$\gamma$ many times. The integral (\ref{cint}) becomes smaller when one
increases the number of turns around $\gamma$. However the estimate
(\ref{ccint}) will hold independently of the number of turns, i.e.
independently of $N$. }
\end{example}

\begin{definition}\label{cdeff3}
Fix a closed 1-form $\omega$ representing $\xi$. The number
\begin{eqnarray}
\ccat(X,\xi)
\end{eqnarray}
is the minimal integer $k$ with the property that there exists $C>0$ such that
for any $N>0$ there exists an open cover $F, F_1, \dots, F_k\subset X$ with the
property that each inclusion $F_j\to X$ is null-homotopic and such that $F$ is
$N$-movable with control $C$ with respect to $\omega$.
\end{definition}

\begin{example}\label{cex2}
{\rm Assume that $X=Y\times S^1$ and $\xi\in H^1(X;\R)$ is such that
$\xi|_{S^1}\not=0$. Then $\ccat(X,\xi)=0$, compare Example \ref{ex1}. }
\end{example}

\begin{lemma}\label{lmless} If $\xi\not=0$ and $X$ is connected then $\ccat(X,\xi)\leq
\cat(X)-1$.
\end{lemma}
This follows from the remark of Example \ref{ex12}
\begin{lemma}\label{zero} If $\xi=0$ then $\ccat(X,\xi)=\cat(X,\xi)=\cat(X)$.
\end{lemma}
The proof of Lemma \ref{zero} is identical to the argument of Example 10.8 in
\cite{farbe3}.

Next we state the homotopy invariance property of $\ccat(X,\xi)$.
\begin{lemma}
Let $\phi: X_1\to X_2$ be a homotopy equivalence, $\xi_2\in H^1(X_2;\R)$, and
$\xi_1=\phi^\ast(\xi_2)\in H^1(X_1;\R)$. Then
\begin{eqnarray}
\ccat(X_1, \xi_1) = \ccat(X_2, \xi_2).
\end{eqnarray}
\end{lemma}
The proof repeats the arguments of Lemma 10.12 in \cite{farbe3}.

It is obvious that in general
\begin{eqnarray}\label{ineqcat}
\ccat(X,\xi) \geq \cat(X,\xi).
\end{eqnarray}
In all examples known to us we have that (\ref{ineqcat}) is an equality.

\begin{remark} {\rm Note that the potentially larger quantity $\ccat(X,\xi)$ can
replace $\cat(X,\xi)$ in the applications to dynamics described in
\cite{farbe2}, \cite{farbe3}, \cite{Lat}, \cite{schu}. This gives an additional
incentive to be interested in the new invariant $\ccat(X,\xi)$. All homotopies
$h_t: F\to X$ which appear in applications to dynamics are induced by gradient
flows of closed 1-forms and hence they satisfy conditions of Definition
\ref{cmovable} with the control constant $C=0$.}
\end{remark}

\begin{example}\label{specific1} {\rm Examining the example described in subsection
\ref{specific} we easily find that for $X=T^2\vee S^1$ one has:
\begin{eqnarray}\label{cases1}
\cat(X,\xi)\, =\, \ccat(X,\xi)\,  = \, \left\{ \begin{array}{ll} 1, &
\mbox{if $\xi\in \R^3 -\ell$},\\
2, & \mbox{if $\xi\in \ell-\{0\}$},\\
3, & \mbox{if $\xi=0$}
\end{array}
\right.
\end{eqnarray}
(as all homotopies described in \S \ref{specific} are in fact with control).
Here $\ell\subset \R^3=H^1(X;\R)$ denotes the set of all classes $\xi\in
H^1(X;\R)$ such that $\xi|_{T^2}=0$.}
\end{example}

\section{Product inequality}

The classical product inequality for the Lusternik - Schnirelmann category
states that
$$\cat(X\times Y)\leq \cat(X)
+ \cat(Y) -1,$$ see \cite{clot}, Theorem 1.37. In this section we show that a
similar inequality holds for $\ccat(X,\xi)$.
\begin{theorem}\label{prodineq}
Let $X$ and $Y$ be finite cell complexes and let $\xi_X\in H^1(X;\R)$ and
 $\xi_Y\in H^1(Y;\R)$ be real cohomology classes. Assume that
\begin{eqnarray}\label{posi}
\ccat(X, \xi_X)>0\quad \mbox{or}\quad \ccat(Y, \xi_Y)>0.
\end{eqnarray}
 Then
\begin{eqnarray}\label{ineqprod}
\ccat(X\times Y, \xi) \leq \ccat(X, \xi_X) + \ccat(Y, \xi_Y) -1,
\end{eqnarray}
where \begin{eqnarray}\xi=\xi_X\times 1 \, +\,  1\times \xi_Y.\end{eqnarray}
\end{theorem}

The $C$-control assumption which appears in Definitions \ref{cmovable} and
\ref{cdeff3} is used in an essential way in  the proof. Theorem \ref{prodineq}
is the main motivation for introducing $\ccat(X,\xi)$ as an alternative to
$\cat(X,\xi)$.

We will start with some auxiliary statements.

\begin{lemma}\label{lmextend} Let $\omega$ be a continuous closed 1-form on a finite cell complex $X$.
Let $F\subset X$ be an open subset which is $N$-movable with respect to
$\omega$ with control $C$ (see Definition \ref{cmovable}). Given a closed
subset $A\subset F$ there exists an open set $F'$ such that $A\subset F'\subset
F$ and a homotopy $H_t: X\to X$, $t\in [0,1]$ satisfying the following
properties: (1) $H_0(x)=x$ for all $x\in X$; (2) for any point $x\in F'$ one
has
\begin{eqnarray}\label{extension}
\int_x^{H_1(x)} \omega \, < \, -N,
\end{eqnarray}
and (3) for any $x\in X$ and for any $t\in [0,1]$ holds
\begin{eqnarray}\label{extension1}
\int_x^{H_t(x)} \omega \, < \, C +1.
\end{eqnarray}
\end{lemma}

The statement of the Lemma can be rephrased as follows: the homotopy which
appears in Definition \ref{movable} can be extended to the whole space $X$ such
that the control condition (i.e. inequality (\ref{ccint})) will hold for all
$x\in X$ and the main inequality (\ref{cint}) will hold on a slightly smaller
set.

\begin{proof}[Proof of Lemma \ref{lmextend}]
Let $h_t: F\to X$, $t\in [0,1]$ be a homotopy given by Definition
\ref{cmovable}. Since $X$ is an ANR and $A\subset X$ is a closed subset, the
homotopy $h_t|_A:A \to X$ can be extended to a homotopy $h'_t: X\to X$, $t\in
[0,1]$, see \cite{Sp}, chapter 1, exercise D2. We have $h'_0(x)=x$ for all
$x\in X$ and $h'_t(x)=h_t(x)$ for all $x\in A$ and $t\in [0,1]$. Using
continuity, we may find a neighborhood $F''\subset F$ of $A$ such that
(\ref{extension}) holds for all $x\in F''$ and such that (\ref{extension1})
holds for all $x\in F''$ and for any $t\in [0,1]$.

Now we want to change the homotopy $h'_t$ such that (\ref{extension1}) holds
for all. Find an open set $F'$ with $A\subset F'\subset \overline{F'}\subset
F''$. Find a continuous function $\phi: X\to [0,1]$ such that $\phi(x)=1$ for
$x\in \overline{F'}$ and $\phi(x)=0$ for $x\in (X-F'')$. Define $h''_t: X\to X$
by $h''_t(x) = h'_{\phi(x)t}(x)$, where $x\in X$ and $t\in [0,1]$. Now it is
quite obvious that the obtained homotopy $h''_t$ satisfies all required
conditions.
\end{proof}

\begin{proof}[Proof of Theorem \ref{prodineq}]
Let $\omega_X$ (correspondingly, $\omega_Y$) be a continuous closed 1-form on
$X$ (correspondingly, $Y$) representing $\xi_X$ (correspondingly $\xi_Y$). Then
$\omega= \omega_X+\omega_Y$ represents $\xi$. More precisely, if $U\subset X$
and $V\subset Y$ are simply connected open subsets and if $\xi_X|_U=df_U$,
$\xi_Y|_V=dg_V$ then $\omega|_{U\times V} = df_U+fd_V$, see \cite{farbe3}, \S
10.2.

Denote $\ccat(X,\xi_X)=r$ and $\ccat(Y, \xi_Y)=s$. There exists $C>0$ such that
for any $N>0$ there exist an open cover $F_X\cup F_1^X\cup \dots\cup F_r^X=X$
with the set $F_X\subset X$ being $(N+C+1)$-movable with respect to $\omega_X$
with control $C$ and with each inclusion $F_j^X\to X$ null-homotopic. Let
$F_Y\cup F_1^Y\cup \dots\cup F_s^Y=Y$ be a cover of $Y$ having similar
properties.

Denote $A_X= X-(F_1^X \cup \dots \cup F_r^X)$. It is a closed subset of $F_X$.
We denote by $H^X_t: X\to X$ the homotopy given by Lemma \ref{lmextend} applied
to $A_X\subset F_X$.

Similarly, consider the set $A_Y= Y-(F_1^Y \cup \dots \cup F_s^Y)$ and denote
by $H^X_t: X\to X$ the homotopy given by Lemma \ref{lmextend} applied to
$A_Y\subset F_Y$.

We obtain that the homotopy $H_t: X\times Y\to X\times Y$ given by
$$H_t(x, y) = (H_t^X(x), H_t^Y(y)), \quad x\in X, y\in Y$$
restricted to a neighborhood $F\subset X\times Y$ of the set $A_X\times Y \cup
X\times A_Y$ satisfies conditions of Definition \ref{cmovable}.

Now the set $F$ together with $F^X_i\times F^Y_j$ where $i=1, \dots r$ and
$j=1, \dots, s$ cover $X\times Y$. It is well known that the union of the sets
$F^X_i\times F^Y_j$ can be covered by $r+s-1$ open sets\footnote{Here we use
our assumption (\ref{posi}).} each null homotopic in $X\times Y$, see
\cite{clot}, \cite{Ja}. Hence we obtain (\ref{ineqprod}).
\end{proof}

\section{Spaces of category zero}

In this section we collect some simple observations about spaces with
$\cat(X,\xi)=0$ for some $\xi\in H^1(X;\R)$. More information can be found in
\cite{schu}, \S3.

\begin{lemma} Let $X$ be a finite CW complex and $\xi\in H^1(X;\R)$. The following properties
are equivalent:

{\rm (1)} $\cat(X, \xi)=0$.

{\rm (2)} $\ccat(X,\xi)=0$.

{\rm (3)} There exists a continuous closed 1-form $\omega$ on $X$ representing
$\xi$ (in the sense of \cite{farbe3}, \S 10.2) and a homotopy $h_t: X\to X$,
where $t\in [0,1]$, such that for any point $x\in X$ one has
\begin{eqnarray}\label{ineq}
\int_x^{h_1(x)}\omega \, <\, 0.
\end{eqnarray}
In (\ref{ineq}) the integral is calculated along the curve $t\mapsto h_t(x)$,
$t\in [0,1]$.

{\rm (4)} For any continuous closed 1-form $\omega$ on $X$ representing $\xi$
there exists a homotopy $h_t: X\to X$, where $t\in [0,1]$, such that for any
point $x\in X$ inequality (\ref{ineq}) holds.
\end{lemma}

\begin{proof} By Definition \ref{deff3}, $\cat(X,\xi)=0$ means that
the whole space $X$ is $N$-movable
 for any $N>0$,
i.e. given $N>0$, there exists a homotopy $H_t: X\to X$, where $t\in [0,1]$,
such that $H_0(x)=x$ and
\begin{eqnarray}\label{large} \int_x^{H_1(x)}\omega <-N
\end{eqnarray}
for any $x\in X$. Hence, (1) implies (3).

 Conversely, given property (3), using compactness of $X$
we find $\epsilon>0$ such that (\ref{ineq}) can be
 replaced by $\int_x^{h_1(x)}\omega <-\epsilon$. Now, one may iterate this
 deformation as follows. The $k$-th iteration is a homotopy $H^k_t: X\to X$,
 where $t\in [0,1]$, defined as follows. Denote by $h_1^{(i)}:X\to X$ the $i$-fold
 composition $h_1^{(i)} = h_1\circ h_1\circ \dots\circ h_1$ ($i$ times). Then for $t\in
 [i/k, (i+1)/k]$ one has
 \begin{eqnarray}
 H^k_t(x) = h_{kt-i}(h_1^{(i)}(x)).
 \end{eqnarray}
 If it is known that $\int_x^{h_1(x)}\omega <-\epsilon$ for any $x\in X$ then
 for the $k$-th iteration one has $\int_x^{H^k_1(x)}\omega <-k\epsilon$
 and (3) follows assuming that $k>N/\epsilon$.
This shows equivalence between (1) and (3).

(4) $\implies$ (3) is obvious. Now suppose that (3) holds for $\omega$ and let
$\omega_1$ be another continuous closed one-form lying in the same cohomology
class, i.e. $\omega_1=\omega +df$ where $f: X\to \R$ is continuous, see
\cite{farbe3}, \S 10.2. Using compactness of $X$ we may find $C$ such that for
any path $\gamma:[0,1]\to X$ one has $|\int_\gamma df|<C$. Fix $N>C$ and apply
equivalence between (1) and (3) to find a homotopy $h_t:X\to X$ with
$\int_x^{h_1(x)}\omega <-N$. Then one has $\int_x^{h_1(x)}\omega_1 <0$, i.e.
(4) holds.

It is obvious that (2) $\implies$ (1). Hence we are left to show that (1)
implies (2). Given (2) fix a deformation as described in (3). Let $C>0$ be such
that for any $x\in X$ and for any $t\in [0,1]$ one has $\int_1^{h_t(x)}\omega
<C$. Then for any iteration $H^k_t: X\to X$ (see above) one has
$\int_x^{H^k_t(x)}\omega <C$ and the result follows.
\end{proof}

In the case when $X$ is a closed smooth manifold a deformation as appearing in
(2) can be constructed as the flow generated by a vector field $v$ on $X$
satisfying $\omega(v)<0$.

The remark of the previous paragraph explains why the following statement can
be viewed as an analogue of the classical Euler - Poincar\'e theorem:

\begin{theorem}\label{ep} $\cat(X,\xi)=0$ implies $\chi(X)=0$.
\end{theorem}
\begin{proof} Suppose that $\cat(X,\xi)=0$. Then $\xi\not=0$, i.e.
the rank $r$ of class $\xi$ is positive, $r>0$. By Lemma \ref{bair} there
exists transcendental bundle $L\in \V_\xi=(\C^\ast)^r$. If $H^q(X;L)\not=0$ for
some $q$ then one may apply Theorem \ref{perfect} with $k=0$ obtaining
$\cat(X,\xi)>0$ and contradicting our hypothesis. Hence $H^q(X;L)=0$ for all
$q$ which implies
$$\chi(X)=\sum_q (-1)^q \dim H^q(X;L)= 0.$$
\end{proof}
Pairs $(X,\xi)$ with $\cat(X,\xi)=0$ form an {\it \lq\lq ideal\rq\rq}\,  in the
following sense:

\begin{lemma}\label{ideal} Let $X_1$ and $X_2$ be finite cell complexes and $\xi_1\in
H^1(X_1;\R)$, $\xi_2\in H^1(X_2;\R)$. If $\cat(X_1, \xi_1)=0$ then
$\cat(X_1\times X_2, \xi)=0$ where $\xi = \xi_1\times 1+ 1\times \xi_2\in
H^1(X_1\times X_2;\R).$
\end{lemma}

\begin{proof} The statement follows directly by applying the definitions.
\end{proof}

\section{An upper bound for $\ccat(M, \xi)$.}\label{sec6}

If $X$ is a connected cell complex and $\xi\in H^1(X;\R)$ is nonzero then
$\cat(X,\xi)\leq \dim X$ and we have seen examples when this inequality is
sharp, i.e. is an equality. However, it can be improved assuming that $X$ is a
manifold.

\begin{theorem}\label{thm3}
Let $M$ be a closed connected smooth $n$-dimensional manifold and let $\xi\in
H^1(M;\R)$ be nonzero. Then\footnote{It is curious to mention that inequality
(\ref{upper}) becomes false if one replaces $\ccat(X,\xi)$ by another
modification of $\cat(X,\xi)$ which was introduced in \cite{FK}, denoted
$\cat^1(X,\xi)$. We will discuss these issues in detail elsewhere.}
\begin{eqnarray}\label{upper} \ccat(M,\xi)\leq n-1.
\end{eqnarray}
\end{theorem}
\begin{proof} Choose a closed 1-form $\omega$ representing $\xi$ having Morse-type zeros and
having no zeros of Morse index 0 or $n$; here we use our assumption that
$\xi\not=0$ and $M$ is connected, see \cite{lev}.

Theorem \ref{thm3} will be proven once we show that for any $N>0$ there exist
closed subsets $A_1, A_2, \dots, A_{n-1}\subset M$ such that each inclusion
$A_j\to M$ is null-homotopic and such that the complement $F=M
-\cup_{j=1}^{n-1}A_j$ is $N$-movable with respect to $\omega$. It is well-known
(see \cite{clot}, Lemma 1.11) that each such $A_j$ can be enlarged to an open
subset $A_j\subset F_j\subset M$ such that $F_j\to M$ is null-homotopic; thus
one obtains an open cover $F\cup F_1\cup \dots\cup F_{n-1}$ satisfying
conditions of Definition \ref{deff3}.

Fix a gradient-like vector field $v$ for $\omega$. As usual, for a zero $p\in
M$ of $\omega$ we denote by $W^s(p)$ and $W^u(p)$ its stable and unstable
manifolds with respect to the flow $x\mapsto x\cdot t$ generated by $-v$. The
closed 1-form $\omega$ \lq\lq locally decreases\rq\rq\ along the flow. More
precisely, this means that $$\int_x^{x\cdot t}\omega <0$$ assuming that $t>0$
and $x$ not a fixed point of the flow; the integral is calculated along the
trajectory of the flow.

We will assume that the stable and unstable manifolds of zeros of $\omega$
intersect transversely; in particular we will require that for any pair of
distinct zeros $p\not=q$ of the same index $\ind (p)=\ind (q)$ one has
$W^s(p)\cap W^u(q)=\emptyset$. Such $v$ exists by the Kupka-Smale theorem.

 Given $N>0$ and a zero $p\in M$ we denote by $A_N(p)$ the set of all points
$x\in W^s(p)$ such that
\begin{eqnarray}
\int_x^p\omega \geq -N.
\end{eqnarray}
As above, the integral is calculated along the integral trajectory $x\cdot t$
where $t\in [0, \infty)$ (or, equivalently, along an arbitrary curve lying in
the stable manifold). Note that $A_N(p)$ is closed in the stable manifold
$W^s(p)$ but it may be not closed in $M$; the closure of $A_N(p)$ in $M$ may
contain zeros of $\omega$ of index greater than $\ind \, p$. Indeed, suppose
that there is an integral trajectory of the flow connecting a zero $q$  with
$p$ where $\ind \, q >\ind \, p$. Then all points of this trajectory except $q$
lie in the stable manifold $W^s(p)$ and the set $A_N(p)$ is not closed for
$\int_q^p \omega \geq -N$.

For each zero $p\in M$ choose a small compact disc $B_N(p)$ containing $p$ in
its interior such that (i) $B_N(p)\cap B_N(p')=\emptyset$ for $p\not= p'$; and
(ii) the set $$C_N(p) = A_N(p)\cup B_N(p)\subset M$$ is contractible.

\begin{center}
\resizebox{6cm}{4.5cm}{\includegraphics[65,289][511, 652]{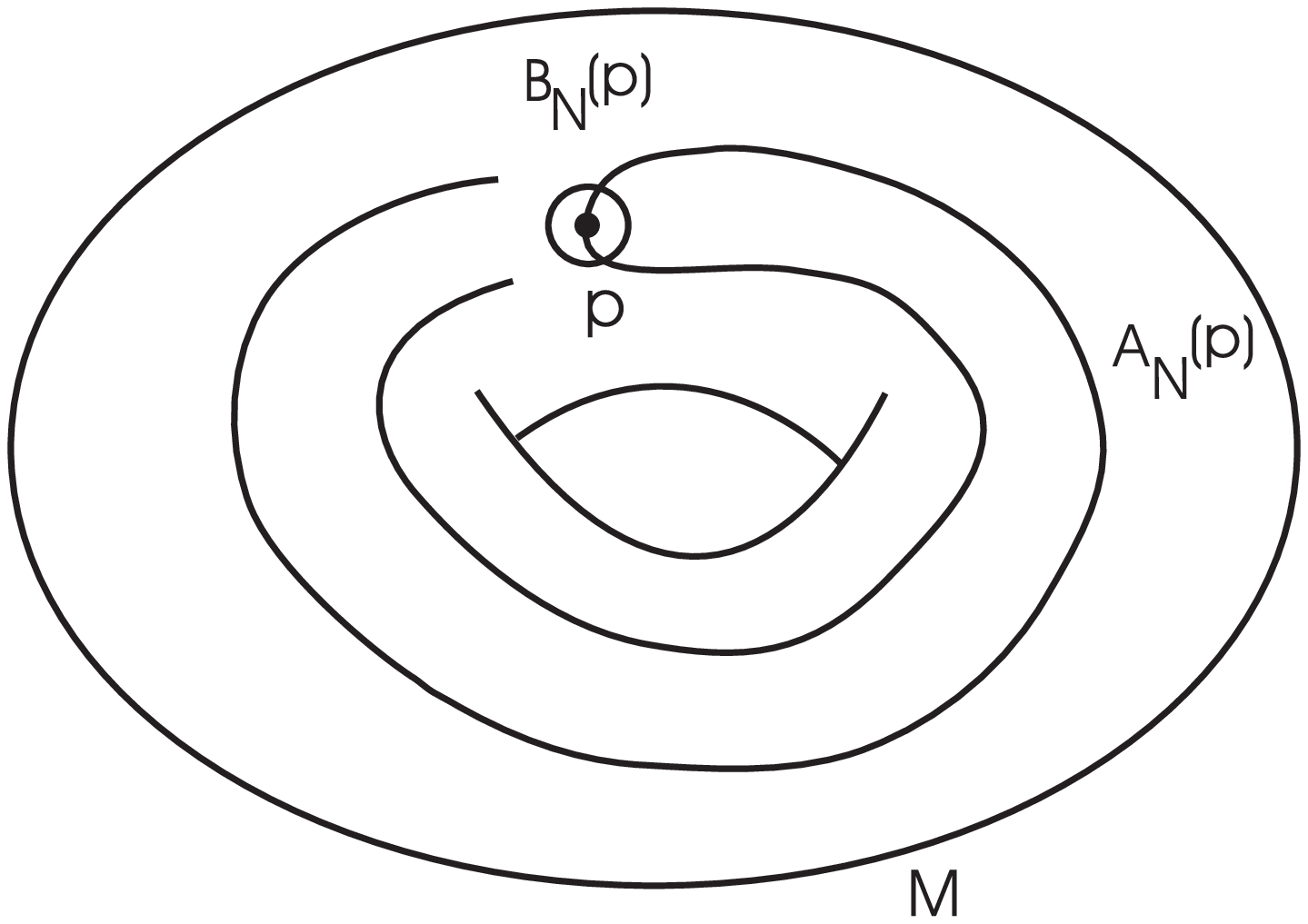}}
\end{center}

Now we claim:{\it

(a) Small discs $B_N(p)$ satisfying (i) and (ii) exist.

(b) The set
$$C'_N(p)= C_N(p) - \bigcup_{\ind (q) >j} \Int B_N(q)$$
is compact. Here $j=\ind p$.}

For $j=1, 2, \dots, n-1$ let $A_j$ be the union of the sets $C'_N(p)$ where $p$
runs over all zeros of $\omega$ of index $j$. We see that each $A_j$ is a
closed subset which is null-homotopic in $M$.

Next we want to show that the complement $$F=M -\bigcup_{\omega_p=0} A_N(p)$$
is $N$-movable with respect to $\omega$. For any point $x\in F$ there exists
$T_x>0$ such that
$$\int_x^{x\cdot T_x}\omega <-N.$$ Hence there is a neighborhood $V_x\subset F$
of $x$ such that
$$\int_y^{y\cdot T_x}\omega <- N$$ for all $y\in V_x$. Without loss of
generality we may assume that for some sequence of points $x_n\in F$ the sets
$V_{x_n}$ form a locally finite cover of $F$. Then we may find continuous
functions $\rho_n: F\to [0,1]$ such that the support of $\rho_n$ is contained
in $V_{x_n}$ and $\max_{n}\rho_n(y) =1$ for any $y\in F$. Then $T=\sum_n
T_{x_n}\rho_n$ is a continuous function on $F$ and for any point $y\in F$ one
has $$\int_y^{yT(y)}\omega < -N.$$ Now a homotopy $h_\tau: F\to M$ as in
Definition \ref{cmovable} (with constant $C=0$) can be defined by the formula
$h_\tau(y) = y\cdot (\tau T(y))$, where $\tau\in [0,1]$.
\end{proof}

\section{Relations with the Bieri-Neumann-Strebel invariant}\label{bieri}

Bieri, Neumann and Strebel introduced in \cite{BNS} a geometric invariant of
discrete groups $G$ which captures information about the finite generation of
kernels of abelian quotients of $G$. In this section we describe a relation
between this invariant and properties of $\cat(X,\xi)$.

Let us recall the definition. We always assume that $G$ is finitely presented
as this is sufficient for our purposes. Let $S(G)$ denote
$(\Hom(G,\R)-\{0\})/\R_+$ where $\R_+$ acts on $\Hom(G,\R)$ by multiplication.
Clearly $S(G)$ is a sphere of dimension $n-1$ where $n$ is the rank of the
abelianization of $G$. Denote by $[\chi]$ the equivalence class of a nonzero
homomorphism $\chi:G\to \R$. The Bieri-Neumann-Strebel invariant associates to
$G$ a subset $\Sigma(G)\subset S(G)$ defined as follows\footnote{We rely on
Theorem 5.1 of \cite{BNS} which states that the definition of $\Sigma(G)$ given
above coincides in the case of finitely presented $G$ with the main definition
of \cite{BNS}.}. Let $X$ be a finite cell complex with $\pi_1(X)=G$ and let $p:
\tilde X \to X$ be the universal abelian cover of $X$. A homomorphism $\chi\in
\Hom(G,\R)$ can be viewed as a cohomology class lying in $H^1(X;\R)$. One has
$\chi\in \Sigma(G)$ if and only if the inclusion $N\to \tilde X$ induces an
epimorphism $\pi_1(N, x_0)\to \pi_1(\tilde X, x_0)$ where $N\subset\tilde X$ is
a connected neighborhood of infinity with respect to $\chi$, see \S
\ref{sec:aj} and Lemma 5.2 from \cite{BNS}.

The following result improves the estimate given by Theorem \ref{thm3}.

\begin{theorem}\label{thmbns}
Let $M$ be a closed connected smooth manifold of dimension $n\geq 5$ and
$G=\pi_1(M)$. If for a nonzero cohomology class $\xi\in H^1(X;\R)$ either
$[\xi]\in \Sigma(G)$ or $[-\xi]\in \Sigma(G)$ then
\begin{eqnarray}\label{ah}
\ccat(M,\xi) &\leq & n-2.
\end{eqnarray}
\end{theorem}
\begin{proof} By Propositions 5.8 and 4.2 of Latour \cite{La} the assumption
$[-\xi]\in \Sigma(G)$ implies that $\xi$ can be realized by a smooth Morse
closed 1-form $\omega$ having no zeros of Morse index $0, 1, n$. Having such
$\omega$ one repeats the argument of the proof of Theorem \ref{thm3} (without
modifications) which leads to (\ref{ah}).

In the case $[\xi]\in \Sigma (G)$ one applies the above argument to $-\xi$ and
obtains a smooth Morse closed 1-form $\omega$ having no zeros of Morse index
$0, n-1, n$. Then one applies the argument of Theorem \ref{thm3}.
\end{proof}

\begin{theorem}\label{nfs}
Let $M$ be a closed connected smooth manifold of dimension $n\geq 5$ such that
$\pi_1(M)$ has no non-abelian free subgroups. Then for any nonzero cohomology
class $\xi\in H^1(X;\R)$ one has
\begin{eqnarray}
\ccat(M,\xi) &\leq & n-2.
\end{eqnarray}
\end{theorem}
\begin{proof}
This follows from the previous statement and from Theorem C of Bieri, Neumann,
Strebel \cite{BNS}. For convenience of the reader we recall that Theorem C of
\cite{BNS} claims that $\Sigma(G)\cup -\Sigma(G) =S(G)$ assuming that $G$ has
no non-abelian free subgroups.
\end{proof}

Here is another statement which is obviously equivalent to Theorem \ref{nfs}:

\begin{theorem}\label{ahh}
Let $M$ be a closed connected smooth manifold of dimension $n\geq 5$ such that
for some nonzero $\xi\in H^1(X;\R)$ one has $\ccat(M,\xi) = n-1.$ Then
$\pi_1(M)$ contains a non-abelian free subgroup.
\end{theorem}

In this theorem one can replace $\ccat(M,\xi)$ by $\cat(M,\xi)$.

Examples of manifolds satisfying $\cat(M,\xi)= \dim M -1$ for some nonzero
$\xi\in  H^1(M;\R)$ are known, see Example 6.5 in \cite{farbe2}. Here $M$ is
obtained from $n$-dimensional torus by adding a handle of index one, i.e
$M=T^n\#(S^1\times S^{n-1})$ and $\xi$ is nonzero restricted to the handle. In
this case it is obvious that $\pi_1(M)$ contains a free subgroup on two
generators.

\begin{theorem}\label{fs}
Let $M$ be a closed connected smooth manifold of dimension $n\geq 5$ and
$G=\pi_1(M)$. If for a nonzero class $\xi\in H^1(X;\R)$ both $[\xi]\in
\Sigma(G)$ and $[-\xi]\in \Sigma(G)$ then
\begin{eqnarray}\label{ah1}
\ccat(X,\xi) &\leq & n-3.
\end{eqnarray}
\end{theorem}
\begin{proof} The proof is similar to the proof of Theorem \ref{thmbns}.
Namely, using Propositions 5.8 and 4.2 of Latour \cite{La} we see that the
assumptions $$[\xi], [-\xi]\in \Sigma(G)$$ imply that $\xi$ can be realized by
a smooth closed 1-form $\omega$ having no zeros of Morse index $0, 1, n-1, n$.
Then one repeats the arguments of the proof of Theorem \ref{thm3} obtaining
(\ref{ah1}).
\end{proof}

Theorem \ref{fs} is similar to Corollary 6.9 from \cite{schu}. It is worth
mentioning that Theorem 6.8 from \cite{schu} states that if $M$ is a closed
smooth manifold of dimension $n\geq 5$ and $\xi\in H^1(M;\R)$ is such that
$[\xi], [-\xi]\in \Sigma(\pi_1(M))$ then $M$ admits a closed 1-form $\omega$
lying in the class $\xi$ which has at most $n-3$ zeros and such that the
gradient flow of $\omega$ has no homoclinic cycles.

\section{Improved upper bound for products}

Combining the product inequality and the upper bound of Theorem \ref{thm3} we
obtain:

\begin{theorem}\label{improved}
Let $M_1, M_2, \dots, M_k$ be closed connected smooth manifolds with $\dim
M_i\geq 2$. Then
\begin{eqnarray}\label{prod}
\ccat(M_1\times M_2\times \dots\times M_k, \xi)\leq 1 \, -2k\,  +\sum_{i=1}^k
\dim M_i
\end{eqnarray}
assuming that the cohomology class $$\xi\in H^1(M_1\times \dots\times M_k;\R)$$
is such that the restriction
$$\xi_i=\xi|_{M_i}\in H^1(M_i;\R)$$
is nonzero for
every $i=1, \dots, k$.
\end{theorem}
\begin{proof} First we consider the case when $\ccat(M_i,\xi_i)=0$ for some
$i=1, \dots, k$. Then the left hand side of (\ref{prod}) vanishes (see Lemma
\ref{ideal}) and inequality (\ref{prod}) is true due to our assumption $\dim
M_i \geq 2$.

If the case when $\ccat(M_i,\xi_i)>0$ for all $i=1, \dots, k$ Theorems
\ref{prodineq} and \ref{thm3} applied repeatedly imply the result by induction.
\end{proof}

\section{Calculation of $\cat(X,\xi)$ for products of surfaces}

\begin{theorem}\label{surface}
Let $M^{2k}$ denote the product $\Sigma_1\times \Sigma_2\times \dots\times
\Sigma_k$ where each $\Sigma_i$ is a closed orientable surface of genus
$g_i>1$. Given a cohomology class $\xi\in H^1(M^{2k};\R)$, one has
\begin{eqnarray}
\cat(M^{2k},\xi)=\ccat(M^{2k},\xi) = 1+ 2r \end{eqnarray} where $r$ is the
number of indices $i\in \{1, 2, \dots, k\}$ such that the cohomology class
$\xi|_{\Sigma_i}\in H^1(\Sigma_i,\R)$ vanishes. In particular
\begin{eqnarray}\label{one1}
\cat(M^{2k},\xi)=\ccat(M^{2k},\xi) = 1
\end{eqnarray}
assuming that $\xi|_{\Sigma_i}\not=0\in H^1(\Sigma_i;\R)$ for any $i=1, \dots,
k$.
\end{theorem}
\begin{proof} If $\Sigma_i$ is a surface of genus $g_i>1$ and $\xi_i\in H^1(
\Sigma_i; \R)$ is nonzero then $\ccat(\Sigma_i,\xi)\leq 1$ by Theorem
\ref{thm3} and $\ccat(\Sigma_i,\xi_i)\geq 1$ by Theorem \ref{ep} (since
$\chi(\Sigma_i) = 2- 2g_i \not=0$). Hence we obtain $\ccat(\Sigma_i, \xi_i)=1.$
Using Theorem \ref{improved} we find that $\ccat(M^{2k}, \xi)=1$ assuming that
$\xi|_{\Sigma_i}\not=0$ for all $i=1, \dots, k$.

Now assume that $\xi_i|_{\Sigma_i}\not= 0$ for $i=1, \dots, k-r$ and
$\xi_i|_{\Sigma_i} =0$ for $i=k-r+1, \dots, k$. Denote
$$M'=\prod_{i=1}^{k-r} \Sigma_i, \quad M''=\prod_{i=k-r+1}^{k} \Sigma_i,
\quad \xi'=\xi|_{M'}, \quad \xi''=\xi|_{M''}=0.$$ As in the previous paragraph
we find $\ccat(M', \xi') =1$. Clearly, $$\ccat(M'',\xi'')=\cat(M'')=2r+1,$$ see
Lemma \ref{zero}. Using the product inequality (Theorem \ref{prodineq}) we
obtain
$$\ccat(M,\xi) \leq \ccat(M', \xi')+\ccat(M'', \xi'') -1= 2r+1.$$

To complete the proof we apply Lemma \ref{clgeq} to get a lower bound. We have
(since the genus of $\Sigma_i$ is greater than $1$)
\begin{eqnarray} \cl(\Sigma_i, \xi_i) =
\left\{
\begin{array} {lll}
0, & \mbox{if} & \xi_i\not=0\in H^1(\Sigma_i;\R),\\
\\
2, & \mbox{if} & \xi_i=0\in H^1(\Sigma_i;\R).
\end{array}
\right.
\end{eqnarray}
By Lemma \ref{clgeq} we obtain $\cl(M^{2k}, \xi) \geq 2r$ where $r$. Hence
$\cat(M^{2k}, \xi) \geq 2r+1$ by Theorem \ref{rephrase}.

This completes the proof.
\end{proof}

\section{Another example}

Let $X$ denote $T^2\vee S^1$ as in subsection \ref{specific}. Formula
(\ref{cases1}) expresses $\ccat(X,\xi)=\cat(X,\xi)$ as function of $\xi\in
H^1(X;\R)=\R^3$.

Consider the product $X^k=X\times \dots \times X$ of $k$-copies of $X$. For any
index $i=1, \dots, k$ denote by $p_i: X\to X^k$ the inclusion $x\mapsto (x_0,
\dots, x_0, x, x_0, \dots, x_0)$ where $x$ stands on the place number $i$ and
on other places is the base point $x_0\in X$. Let $q_i: T^2\to X^k$ be the
composition of the inclusion $T^2\to X$ and of $p_i: X\to X^k$.

\begin{theorem}\label{last} For any $\xi\in H^1(X^k;\R)$ one has
\begin{eqnarray}
\cat(X^k,\xi)=\ccat (X^k,\xi) = 1+a(\xi)+2b(\xi)
\end{eqnarray}
where $a(\xi)$ denotes the number of indexes $i\in \{1, \dots, k\}$ such that
$p_i^\ast(\xi)\not=0$ and $q_i^\ast(\xi)=0$, and $b(\xi)$ denotes the the
number of indexes $i\in \{1, \dots, k\}$ such that $p_i^\ast(\xi)=0$.
\end{theorem}

\begin{proof} Denote by $X_i$ the $i$-th factor $X$ in the product $X^k=X\times
\dots \times X$. Let $\xi_i\in H^1(X_i;\R)$ denote $p_i^\ast(\xi)$. Finally we
denote by $\ell_i\subset H^1(X_i;\R)$ the set of cohomology classes such that
their restriction onto the torus $T^2\subset X_i$ vanishes. We have
\begin{eqnarray}\label{cases2}
\cat(X_i,\xi_i)\, =\, \ccat(X_i,\xi_i)\,  = \, \left\{ \begin{array}{ll} 1, &
\mbox{if $\xi_i\notin \ell_i$},\\
2, & \mbox{if $\xi_i\in \ell_i-\{0\}$},\\
3, & \mbox{if $\xi_i=0$},
\end{array}
\right.
\end{eqnarray}
see (\ref{cases1}). It is easy to check directly using the definitions that
\begin{eqnarray}
\cl(X_i,\xi_i)\, = \, \left\{ \begin{array}{ll} 0, & \mbox{if $\xi_i\notin \ell_i$},\\
1, & \mbox{if $\xi_i\in \ell_i-\{0\}$},\\
2, & \mbox{if $\xi_i=0$}.
\end{array}
\right.
\end{eqnarray}
Applying Theorem \ref{prodineq} we find
$$\cat(X^{2k},\xi) \leq \ccat(X^{2k},\xi) \leq 1+a(\xi)+2b(\xi).$$
Applying Theorem \ref{rephrase} and Lemma \ref{clgeq} inductively we obtain the
inverse inequality
$$\cat(X^{2k},\xi) \geq 1 +a(\xi) +2b(\xi).$$
This completes the proof.
\end{proof}

\section{Questions}

Finally we raise several challenging questions which are inspired by problems
discussed in this article.

\begin{question} Is it always true that $\ccat(X,\xi)=\cat(X,\xi)$?
\end{question}

In all examples discussed in this paper the function $\xi\mapsto \cat(X,\xi)$
was upper semi-continuous, i.e. the sets of the form $\{\xi\in H^1(X;\R); \,
\cat(X,\xi)\geq r\}$ were closed.

\begin{question} Is it true in general that the functions $H^1(X;\R) \to \R$ given by
$$\xi\mapsto \cat(X,\xi),\quad\xi\mapsto \ccat(X,\xi)$$ are upper
semi-continuous?
\end{question}

\bibliographystyle{amsalpha}

\begin{thebibliography}{99}

\bibitem{BG} R. Bieri, R. Geoghegan, \textit{Connectivity properties of group
actions on non-positively curved spaces}, Mem. Amer. Math. Soc. ~161(2003), no.
765.

\bibitem{BNS} R. Bieri, W. Neumann, R. Strebel, \textit{A geometric invariant
of discrete groups}, Invent. Math. 90(1987), 451 - 477.

\bibitem{B} N. Bourbaki, \textit{Commutative algebra}, Hermann and Addison -
Wesley Publishing Company, 1972


\bibitem{clot} O. Cornea, G. Lupton, J. Oprea, D. Tanr\'e, \textit{Lusternik -
Schnirelmann category}, AMS, 2003.


\bibitem{farbe2}M. Farber, \textit{Zeros of closed 1-forms, homoclinic orbits
and Lusternik-Schnirelman theory},  Topol.~Methods Nonlinear Anal.~19 (2002),
123-152.

\bibitem{farbe4} M. Farber, \textit{Lusternik -- Schnirelman theory and dynamics,}
"Lusternik - Schnirelmann
Category and Related Topics", Contemporary Mathematics, ~ 316(2002), 95 - 111.


\bibitem{FK} M. Farber, T. Kappeler, \textit{
Lusternik - Schnirelman theory and dynamics, II}, Proceedings of the Steklov
Institute of Mathematics, vol. 247(2004), issue 4, 232 - 245.

\bibitem{farbe3} M. Farber, \textit{Topology of closed one-forms}. Mathematical
Surveys and Monographs, 108. American Mathematical Society, Providence, RI,
2004.

\bibitem{FS} M. Farber, D. Sch\"utz, \textit{Moving homology classes to
infinity}, preprint, 2005, to appear in {\it Forum Mathematicum}.

\bibitem{FS1} M. Farber, D. Sch\"utz, \textit{Closed 1-forms with at most one zero},
preprint, 2005, to appear in {\it Topology}.


\bibitem{Ja} I.M. James, \textit{On the category in the sense of Lusternik-Schnirelman},
Topology, \textbf{17}(1978), pp. 331 - 348.

\bibitem{La} F. Latour, \textsl{Existence de 1-formes ferm\'ees non
singuli\`eres dans une classe de cohomologie de de Rham}, Publ. IHES 80 (1994),
135-194.

\bibitem{Lat} J. Latschev, {\em Flows with Lyapunov 1-forms and a
generalization of Farber's Theorem on homoclinic cycles}, International
Mathematical Research Notices, 5(2004), 239-247.


\bibitem{lev} G. Levitt, \textit{1-formes ferm\'ees singuli\`eres et
groupe fondamental}, Invent. Math., \textbf{88}(1987), 635 - 667.


\bibitem{neumann} W. Neumann, \textit{Signature related invariants of
manifolds. I. Monodromy and $\gamma$-invariants}, Topology, {\bf{18}}(1979),
147-172.

\bibitem{noviko} S. Novikov, \textit{Multi-valued functions and functionals. An
analogue of Morse theory.} Soviet. Math. Doklady, ~24(1981), 222-226.

\bibitem{schu} D. Sch\"utz, \textit{On the Lusternik-Schnirelman theory of
a real cohomology class}, Manuscripta math., ~113(2004), 85-106.



\bibitem{Sp} E. Spanier, \textit{Algebraic topology}, McGraw-Hill Book Company,
1966



\end{thebibliography}

\end{document}